\documentclass{article}

\usepackage{authblk}
\usepackage{amsthm}
\usepackage{amssymb,amsmath,amsfonts,microtype,mathrsfs}
\usepackage{indentfirst}
\usepackage{graphicx}
\usepackage{pdfsync}
\usepackage{comment}
\usepackage{enumitem}
\usepackage[colorlinks=true,         
            linkcolor=blue,          
            urlcolor=blue,           
            citecolor=blue]{hyperref}

\usepackage{amssymb}
\usepackage{hyperref}

\usepackage{microtype}
\usepackage{esint}

\usepackage[colorlinks=true,         
            linkcolor=blue,          
            urlcolor=blue,           
            citecolor=blue]{hyperref}
\theoremstyle{definition} 
\newtheorem{thm}{Theorem}
\newtheorem{cor}{Corollary}[section]
\newtheorem{lem}{Lemma}[section]
\newtheorem{prop}{Proposition}[section]

\newtheorem*{conj}{Conjecture}

\newtheorem*{thmB}{Theorem B}
\newtheorem*{thmA}{Theorem A}
\newtheorem*{thmC}{Theorem C}
\newtheorem*{thmD}{Theorem D}
\newtheorem*{thmE}{Theorem E}
\newtheorem*{thmF}{Theorem F}
\newtheorem*{thmB'}{Theorem B$^\prime$}
\newtheorem*{thmC'}{Theorem C$^\prime$}

\theoremstyle{definition}

\newtheorem*{rem}{Remark}
\theoremstyle{remark}

\newcommand{\R}{\mathbb{R}}

\newcommand{\Z}{\mathbb{Z}}
\newcommand{\N}{\mathbb{N}}

\newcommand{\C}{\mathbb{C}}
\newcommand{\F}{\mathbb{F}_q[t]}

\newcommand{\tr}{\text{tr}}
\renewcommand\phi\varphi

\newcommand{\1}{\mathbf{1}}

\numberwithin{equation}{section}

\title{Sárközy's theorem in $\mathbb{F}_q[t]$ via the van der Corput property}

\author{Steve Fan\thanks{Steve.Fan@uga.edu}\quad and Andrew Lott\thanks{Andrew.Lott@uga.edu}}
\affil{Department of Mathematics\\ University of Georgia\\ Athens, GA 30602}

\date{}

\begin{document}

\maketitle

\begin{abstract}
Fix a positive prime power $q$, and  let $\F$ be the ring of polynomials over the finite field $\mathbb{F}_q$ with $\text{char}(\mathbb{F}_q)>2$. Suppose $A \subseteq \{f \in \F : \deg f \leq N\}$ contains no pair of elements whose difference is of the form $P-1$ with $P$ irreducible. Adapting Green’s approach to Sárközy’s theorem for shifted primes in $\Z$ using the van der Corput property, we show that
\[
|A| \ll q^{(N+1)(11/12+o(1))},
\]
improving upon the bound $O\big(q^{(1-c/\log N)(N+1)}\big)$ due to L\^{e} and Spencer.  An important distinction between Green's argument and ours lies in the properties of exponential sums over function fields, which differ in several interesting ways from their number-field counterparts.

\end{abstract}
\tableofcontents

\section{Introduction}

For $N\in \mathbb{N}$, let $D(N)$ denote the size of the largest subset of $\{1,...,N\}$ containing no two elements differing by $p-1$ with $p$ prime. The study of $D(N)$ can be traced back to a question of Erdős.\footnote{This is according to Sárközy \cite{SarkozyIII1978}, but he provides no reference.}

\begin{conj}[Erdős] We have
\[\lim_{N\to\infty}\frac{D(N)}{N}=0.\]
\end{conj}

In $1978$, Sárközy \cite{SarkozyIII1978} confirmed this conjecture, proving the first quantitative upper bound on $D(N)$ by utilizing the density increment strategy of Roth \cite{SarkozyIII1978}. 
\begin{thmA}[Sárközy, 1978] For sufficiently large $N$,
\[D(N)\ll\frac{N}{(\log\log N)^{2+o(1)}}.\]
\end{thmA}

Thirty years later, Lucier \cite{Lucier08} provided the first improvement on S\'{a}rkozy's bound, showing that 
\[D(N)\ll \frac{N}{(\log\log N)^{(1+o(1))\log\log\log\log N}}\]
for sufficiently large $N$. Shortly afterwards, Ruzsa and Sanders \cite{RuzsaSanders2008} established a quasipolynomial upper bound on $D(N)$ by using an energy increment strategy. 
\begin{thmB}[Ruzsa--Sanders, 2008]
There exists a constant $c>0$ such that 
\[D(N)\ll Ne^{-c(\log N)^{1\slash 4}}.\]
\end{thmB}

In 2020, Wang \cite{Wang2020} was able to optimize the approach of Ruzsa--Sanders to achieve an upper bound of the form $Ne^{-c(\log N)^{1\slash 3}}$. Even assuming GRH, the proof of Ruzsa--Sanders can only yield a bound of quasipolynomial shape.

In 2023, Green established a power-saving bound on $D(N)$ by proving that the shifted primes $p-1$ satisfy a strong form of the van der Corput property, which is strikingly different from the more classical increment methods \cite{Green23}.

\begin{thmC}[Green, 2023]\label{GreenUnconditional}
There exists a constant $c>0$ such that 
\[D(N)\ll N^{1-c}.\]
\end{thmC}

Combining Green's arguments with their explicit version of Bombieri’s refinement of Gallagher’s log-free large sieve density estimate for Dirichlet $L$-functions near the line $\sigma=1$, Thorner and Zaman \cite{ThornerZaman2024} recently obtained unconditionally the admissible constant $c=10^{-18}$ for Theorem \ref{GreenUnconditional}.

Under GRH Green’s method yields an improved exponent, and the proof becomes significantly easier \cite{Green}.

\begin{thmD}[Green, 2023]\label{GreenConditional}
Assume GRH. For any $\epsilon>0$, we have 
\[
D(N) \ll_{\epsilon} N^{11/12+\epsilon}.
\]
\end{thmD}

The state-of-the-art lower bound for $D(N)$ is due to a clever construction of Ruzsa from $1984$ \cite{Ruzsa1984}. 

\begin{thmE}[Ruzsa, 1984]\label{lowerbound}
For sufficiently large $N$, 
\[D(N)\geq \exp\left(\big(\log2+o(1)\big)\frac{\log N}{\log \log N}\right).\]
\end{thmE}

Green \cite{Green23} expects that this is close to the true lower bound for $D(N)$, which would mean that the upper bound in Theorem C is still much larger than it needs to be.

In the present paper we establish analogues of Theorem C and Theorem E in $\F$, the ring of polynomials of the finite field $\mathbb{F}_q$. In this setting $\F$ plays the role of the integers, irreducible polynomials play the role of the primes, and $\{f\in \F: \deg f\leq N\}$ plays the role of $\{1,...,N\}$. In particular, for $N\in \N$, let $D_q(N)$ denote the size of the largest subset of $\{f\in \F: \deg f\leq N\}$ containing no two elements differing by $P-1$ with $P$ irreducible.

\begin{thm}\label{MainTheorem}
Let $\epsilon>0$ and let $N$ be sufficiently large. Then, 
\[ D_q(N) \ll_{q,\epsilon} q^{(N+1)(11/12+\epsilon)},
\]
provided $\text{char}(\mathbb{F}_q)>2$. On the other hand, we have
\[D_q(N)\ge\exp\left((\gamma_q\log2+o_q(1))\frac{N\log q}{q\log N}\right),\]
where 
\[\gamma_q:=
\begin{cases}
1/2&\text{ if $\text{char}(\mathbb{F}_q)>2$},\\
1&\text{ if $\text{char}(\mathbb{F}_q)=2$}.
\end{cases}\]
\end{thm}

An important difference between Green's argument and ours lies in the nuances of exponential sum estimates over function fields, which are more rigid than their number field counterparts. This same rigidity also means that Green's smoothing kernel is not needed in our argument. Using character sum estimates \cite{KeatingRudnick2014}, we are also able to derive an ``all degrees" version of the exponential sum estimate of Hayes \cite{Hayes66}, which plays a key role in our proof. 

We remark that Theorem \ref{MainTheorem} still holds with $P-1$ replaced by $P-a$ for any given $a\in\mathbb{F}_q^{\times}$. The upper bound supplied by Theorem \ref{MainTheorem} improves on the bound $O\big(q^{(1-c/\log N)(N+1)}\big)$ with some constant $c>0$ obtained by L\^{e} and Spencer \cite{LeSpencer2011}. We achieved the upper bound by applying Weil's analogue of GRH and adapting the proof of Theorem D. The lower bound in Theorem \ref{MainTheorem} is essentially the same as that of L\^{e} and Spencer but with a better asymptotic constant in the $q$ aspect. It arises from a straightforward adaptation of Ruzsa's classical construction. Since L\^{e} and Spencer omitted the proof of their lower bound, we include a proof of ours for the sake of completeness.

On the other hand, substantial progress has been made on the analog of the Furstenberg--Sárközy theorem for differences of $k$th powers $g^k$ in $\F$, where $k\ge2$ is an integer. Indeed, Green \cite{GreenSquares2017} proved that if $A\subseteq\{f\in \F: \deg(f)\leq N\}$ contains no two distinct elements differing by $g^k$ for any $g\in \F$, then 
\[|A|\le 2q^{(N+1)(1-c_{q,k})}\]
with some explicit constant $c_{q,k}>0$ depending on $q,k$. In light of a construction of Ruzsa \cite{Ruzsa1984Squares} (adapted for $\F$ by Link \cite{Link2008}), the shape of this bound is essentially the best one could achieve. Green's proof of this remarkable result borrows ideas from the seminal work of Croot--Lev--Pach \cite{CrootLevPach2017}.

\begin{rem}
The astute reader may have observed the superior lower bound for $D_q(N)$ in our main result in the case where $\mathbb{F}_q$ has characteristic 2 and the exclusion of this case in the statement of the upper bound. As in many results in function fields, the case $\text{char}(\mathbb{F}_q)=2$ needs some special care. Shortly after finishing this paper, the authors learned from A. Kowalska \cite{Kowalska26} that she had independently obtained a result on Sárközy's theorem in $\mathbb{F}_2[t]$. Most notably, she obtained the better exponent $7/8$ for $D_2(N)$. We expect her method to also yield improvements in $\mathbb{F}_{2^{\ell}}[t]$ with $\ell$ a fixed positive integer.
\end{rem}

\section{Notation}
Throughout the paper, we write $\mathbb{F}_q$ for the finite field of order $q=p^{\ell}$ and characteristic $\text{char}(\mathbb{F}_q)=p$, with $p$ prime and $\ell$ a positive integer, and denote by $\F$ the ring of polynomials over $\mathbb{F}_q$. The letters $k,n,N$ will always represents positive integers, the letters $f,g,h$ will denote polynomials in $\F$, and the letter $P$ will be reserved for irreducible polynomials. For any $X\in\R$, we write $\lfloor X\rfloor$ for the integer part of $X$, which is defined as the greatest integer not exceeding $X$. Given $f\in\F$, we denote by $(f)$ the principal ideal generated by $f$. For any $f,g\in \F$, we write $\gcd(f,g)$, or simply $(f,g)$, for the greatest, monic common divisor of $f$ and $g$, and $[f,g]$ for the least, monic common multiple of $f$ and $g$.

We adopt the notational conventions of Liu and Wooley \cite{LiuWooley2010} for our function field setup. Let $K$ be the field of fractions of $\F$. For every $r=f\slash g \in K$ with $f,g\in\F$, we define the {\it order} of $r$ by $\text{ord}(r):=\deg f-\deg g$ and the {\it norm} $\langle \cdot \rangle: K\to [0,\infty)$ by 
\[\langle r \rangle :=q^{\text{ord}(r)}.\]
Here, we follow the convention that $\deg 0=-\infty$ so that $\langle 0 \rangle = q^{-\infty}=0$. It is not hard to check that $\langle \cdot \rangle$ is indeed a norm, and one obtains the completion $K_\infty$ of $K$ with respect to this norm. Both the order and the norm can be extended to $K_\infty$: given any $x\in K_{\infty}$ expressed as a Laurent expansion of the form
\[x=\sum_{n\in\Z}a_nt^n,\]
where each $a_n\in\mathbb{F}_q$ and $a_n=0$ for all but finitely many $n\in\N$, we define $\text{ord}(x):=\max\{n\in \Z:a_n\neq 0\}$ and $\langle x \rangle: = q^{\text{ord}(x)}$. And as in complex analysis, the {\it residue} $\text{res}(x)$ of $x$ is defined by $\text{res}(x):=a_{-1}$.

Next, we introduce the analogues for $\R\slash \Z$ and $\mathbb{Q}\slash \Z$:
\[\mathbb{T}:=K_\infty \slash \F=\{x\in K_\infty: \langle x \rangle <1\}\]
and 
\[K\slash \F:=\mathbb{T}\cap K.\]
To make the analogy with the number field case clearer, for every $R\in \R$ we also introduce the notation 
\[\widehat{R}:=q^R.\]

In order to do Fourier analysis on $\F$, we require an analogue of the classical complex exponential function $e^{2\pi i x}$ that works in the function field setting. To this end, we define, for any $x\in K_\infty$, 
\[e(x)=\exp\left(2\pi i\frac{\text{tr}(\text{res}(x))}{p}\right),\]
where we recall $q=p^\ell$, and where the field trace is given by $\tr(a):=a+a^p+\cdots +a^{p^{\ell-1}}$. The function $e\colon K_\infty\to\C^{\times}$ defined this way is easily seen to be an additive character enjoying many similar properties as the classical complex exponential function, with the most important being the orthogonality relation (see Lemma \ref{orthogonalityrelation} below), which forms the basis for Fourier analysis on $\F$.

We will also adopt the standard Landau--Vinogradov asymptotic notation such as $O$, $o$, $\ll$ and $\gg$ from analytic number theory. Given complex-valued functions $F,G$ of the variable $N$ in a certain range, we write $F=O(G)$ or $F\ll G$ if there is a constant $C>0$ such that $|F|\le C|G|$ for all $N$ in the considered range. In our context, the constant $C$ is allowed to depend on $q$, and it may also depend on other parameters occasionally. The relation $F\gg G$ is then used interchangeably with $G=O(F)$, and the relation $F=o(G)$ is interpreted as $F/G\to0$ as $N\to\infty$. Lastly, for any subset $A\subseteq\F$, the indicator function $1_{A}$ of $A$ is defined by $1_{A}(f)=1$ if $f\in A$ and $1_{A}(f)=0$ otherwise. Analogously, we define, for any logical statement $S$, $1_{S}=1$ if $S$ is true and $1_{S}=0$ otherwise.

In addition, we will need analogues of a few classical arithmetic functions on $\Z$. For any $f\in \F\setminus\{0\}$ with irreducible factorization $f=uP_1^{\alpha_1}\cdots P_m^{\alpha_m}$, where $P_1,...,P_m\in\F$ are distinct monic, irreducible polynomials, $\alpha_1,...,\alpha_m\in\N$, and $u\in\mathbb{F}_q^{\times}$, we define the divisor function $\tau$ by 
\[\tau(f):=\prod_{i=1}^m (\alpha_i+1),\]
the Euler totient function $\phi$ by
\[\phi(f):=\#(\F/(f))^{\times}=\langle f \rangle\prod_{i=1}^{m}\left(1-\frac{1}{\langle P_i \rangle}\right),\]
and the M\"{o}bius function $\mu$ by
\[\mu(f):=
\begin{cases}
(-1)^m & \text{ if $\alpha_1=\cdots=\alpha_m=1$},\\
0 & \text{ otherwise}.
\end{cases}
\]
We adopt the convention that $\tau(f)=\phi(f)=\mu(f)=1$ for any $f\in\mathbb{F}_q^{\times}$. For the same $f$ described above, we also introduce the von Mangoldt function 
\[\Lambda(f):=
\begin{cases}
\deg P_1&\text{ if $m=1$},\\
0&\text{ otherwise},
\end{cases}\]
along with its close cousin
\[\Lambda'(f)=
\begin{cases}
\deg f&\text{ if } f\text{ is irreducible},\\
0&\text{ otherwise}.
\end{cases}\]
When it comes to counting irreducible polynomials, let $\pi(n)$ count the number of monic, irreducible polynomials of degree $n\ge1$. Given a non-constant $f\in\F$ and $a\in\F$, let 
\[\pi(n;f,a):=\#\{P\in\F\text{~monic, irreducible of degree~}n\colon P\equiv a\hspace*{-2.5mm}\pmod{f}\}.\]
In particular, $\pi(n)=\pi(n;1,1)$. In the case $\gcd(a,f)=1$, the proof of \cite[Theorem 4.8]{Rosen02} (with uniformity in both $f$ and $a$ taken into account) shows that 
\begin{equation}\label{eq:PNT}
\pi(n;f,a)=\frac{1}{\phi(f)}\cdot\frac{q^n}{n}+O\left(\frac{q^{n/2}}{n}\deg f\right),
\end{equation}
which is a consequence of the Generalized Riemann Hypothesis (GRH) for $\F$ proved by A. Weil in 1948.

Finally, given non-constant $f\in\F$, a {\it Dirichlet character} $\chi\pmod{f}$ is a completely multiplicative function from $\F$ to $\{z\in\C\colon |z|=1\}\cup\{0\}$, with the properties that $\chi(a)=0$ if and only if $\gcd(a,f)>1$ and that $\chi(a)=\chi(b)$ for any $a,b\in\F$ with $a\equiv b\pmod{f}$. The trivial character $\chi_0\pmod{f}$ is the character satisfying $\chi_0(a)=1$ for all $a\in\F$ with $\gcd(a,f)=1$. There are exactly $\phi(f)$ Dirichlet characters modulo $f$, and they satisfy the following well-known orthogonality relations:
\[\frac{1}{\phi(f)}\sum_{\chi\hspace*{-2.5mm}\pmod{f}}\bar{\chi}(a)\chi(b)=1_{a\equiv b\hspace*{-2.5mm}\pmod{f}}=
\begin{cases}
1&\text{ if $a\equiv b\hspace*{-2.5mm}\pmod{f}$},\\
0&\text{ otherwise},
\end{cases}\]
for any $a,b\in\F$ coprime to $f$, and
\[\frac{1}{\phi(f)}\sum_{a\in\F/(f)}\bar{\chi}_1(a)\chi_2(a)=1_{\chi_1=\chi_2}=
\begin{cases}
1&\text{ if $\chi_1=\chi_2$},\\
0&\text{ otherwise},
\end{cases}\]
for any characters $\chi_1,\chi_2\pmod{f}$. Moreover, GRH implies 
\begin{equation}\label{eq:charactersum}
|\psi(n,\chi)|\le(\deg f-1)q^{n/2}+\sum_{\substack{k\mid n\\ k\le n/2}}\pi(k)k=\left(\deg f+O\left(d(n)\right)\right)q^{n/2},    
\end{equation}
for all $\chi\ne\chi_0$ \cite{KeatingRudnick2014}, where 
\begin{equation}\label{eq:psi(n,chi)}
\psi(n,\chi):=\sum_{\substack{\deg g=n\\g~\text{monic}}}\Lambda'(g)\chi(g),   
\end{equation}
and $d(n)$ counts the number of positive divisors of $n$. 

\section{Fourier analysis in function fields}

In \cite{Green}, Green applies the circle method to handle certain difficult exponential sums over the primes. In analogy with his approach, we will require a suitable analogue for exponential sums over irreducibles in $\F$. In this section we introduce several standard results on Fourier analysis and character sums in function fields whose proofs can be found in various sources  \cite{Car96,EffingerHayes91,Hayes66,Kubota74,LiuWooley2010}.

We start with the following orthogonality relation for exponential functions over $\F\slash (h)$. 
\begin{lem}[Orthogonality relation]\label{orthogonalityrelation}
Let $h\in \F$ be nonzero and $f\in \F/(h)$. Then,
\[
\sum_{g\in \F/(h)}e\left(\frac{fg}{h}\right)
=
\begin{cases}
0&\text{ if } f\not\equiv 0\pmod{h},\\
\langle h \rangle&\text{ if }f\equiv 0\pmod{h}.
\end{cases}
\]
\end{lem}

Let $h\in \F$ be nonzero. We define the Fourier transform of a function $F: \F\slash (h)\to \C$ as follows: 
\[\widehat{F}(g)=\frac{1}{\langle h \rangle }\sum_{r\in \F\slash (h)}F(r)e(-rg \slash h).\]
Then, the orthogonality relation immediately yields the Fourier inversion formula and Plancherel's identity. 

\begin{lem}[Fourier inversion formula]
Let $h\in \F$ be nonzero, and let $F: \F\slash (h)\to \C$. Then, 
\[F(g)=\sum_{r\in \F\slash (h)}\widehat{F}(r)e(rg\slash h).\]
\end{lem}
\begin{lem}[Plancherel's identity]
Let $h\in \F$ be nonzero, and let $F: \F\slash (h)\to \C$. Then, 
\[\sum_{r\in \F\slash (h)}|F(r)|^2=\langle h \rangle \sum_{r\in \F\slash (h)}|\widehat{F}(r)|^2.\]
\end{lem}
With this orthogonality relation in hand, we can also deal with \textit{Ramanujan sums} over $\F$. For $s\in \F$ and $f\in\F$, define the Ramanujan sum
\[c_s(f):=\sum_{a\in (\F\slash (s))^\times}e(af\slash s).\] 
 
The following lemma on $c_s(f)$ is well-known.
\begin{lem}\label{RamanujanSumFormula}
Let $s\in \F$ and $f\in \F\slash (s)$. Then $ c_s(f)$ is multiplicative in $s$. Moreover, we have 
\[c_s(f)=\mu\left(\frac{s}{(s,f)}\right)\frac{\phi(s)}{\phi\big(\frac{s}{(s,f)}\big)}.\]
\end{lem}

The following corollary is immediate. 

\begin{cor}\label{ramanujansumbound}
Let $s\in \F$ and $f\in \F\slash (s)$. Then, 
\[|c_s(f)|\leq \langle (s,f)\rangle.\]
\end{cor}

Next up, we define the Gauss sum associated to a Dirichlet character $\chi\pmod{r}$ by
\[G(\chi):=\sum_{a\in \F\slash (r)}\chi(a)e\left(\frac{a}{r}\right).\]
In particular, we have $G(\chi_0)=c_r(1)=\mu(r)$. We will need the following simple upper bound for $G(\chi)$ when $\chi$ is nontrivial. 

\begin{lem}\label{lem:Gausssum}
For any nontrivial character $\chi\pmod{r}$, we have
\[|G(\chi)|\le\langle r\rangle^{1/2}.\]
\end{lem}

We will also make heavy use of the function field analogue for Dirichlet's approximation theorem, which can be found in \cite{Pollack13}. 
 
\begin{lem}[Dirichlet]
Let $x\in \mathbb{T}$ and $R>0$. Then there exist $b, r\in \F$ coprime which satisfy $\langle r \rangle \leq \widehat{R}$ and
\[\left\langle\frac{b}{r}-x \right\rangle<(\langle r \rangle \widehat{R})^{-1}.\]
\end{lem}

Whenever doing the circle method, one runs into the following fundamental bound, which comes from the geometric sum formula: 
\[\left|\sum_{n=0}^Ne^{2\pi i nx}\right|\ll\min\{\|x\|^{-1},N\},\]
where $\|x\|$ denotes the distance of $x\in\R$ to the nearest integer. The function field analogue of this turns out to be even nicer, in a somewhat surprising way.
\begin{lem}\label{exponentialsumformula} Let $x\in \mathbb{T}$ and $N\in \mathbb{N}$. Then, 
\[\sum_{\deg(f)\leq N }e(fx)=
\begin{cases}
q^{N+1} &\text{ if } \text{ord } (x)  <{-N-1},\\
0&\text{ otherwise}.
\end{cases}
\]
\end{lem}

We also have the following corollary, which will be used in Section \ref{sec:proof of Prop 5.1}. 

\begin{cor}\label{exponentialsumformulaV2} Let $x\in \mathbb{T}$ and $N\in \mathbb{N}$. Then, 
\[\sum_{\substack{\deg(f)=N\\f~\text{monic}} }e(fx)=
\begin{cases}
q^{N}e(t^Nx) &\text{ if } \text{ord } (x)  <{-N},\\
0&\text{ otherwise}.
\end{cases}
\]
\end{cor}
\section{The van der Corput property in function fields}
The rest of the paper will be devoted to proving that the shifted irreducible elements of $\F$ satisfy the van der Corput property with a power saving bound. This will almost immediately imply the main theorem.  
\begin{prop}[The van der Corput property]\label{vdc}
Let $\epsilon>0$ and $N\in \N$. There exist real coefficients $a_0$ and $\{a_{P-1}: P\in \F\}$ such that the cosine polynomial $T: K_\infty \to \R$ given by 
\[T(x)=a_0+\sum_{
\substack{\langle P\rangle \leq \widehat{N}}}a_{P-1}\Re\left[e\left(x(P-1)\right)\right] \]
satisfies the following three properties:

\begin{center}
\begin{tabular}{rl}
(1) & $T(0) = 1$, \\
(2) & $T(x) \geq 0$ for every $x \in K_\infty$, \\
(3) & $0<a_0 \ll_{\epsilon,q} \widehat{N}^{-1/12+\epsilon}$. \\
\end{tabular}
\end{center}

\end{prop}

We introduce one more piece of notation. Let $T$ be defined as in Proposition \ref{vdc}. Define 
\[\widetilde{T}(x)=a_0+\sum_{
\substack{\langle P\rangle \leq \widehat{N}}}a_{P-1}e\left(x(P-1)\right). \]
We think of $\widetilde{T}$ as the Fourier series of a function supported on $0$ and the shifted irreducibles in $\F$. We are then going to consider $A\subseteq\{f\in \F: \deg f \leq N\}$ and the expression 
\[\sum_{g\in \F/ (t^{N+1})}\Bigg|\sum_{\langle f \rangle \leq \widehat{N}}\1_A(f)e\left(fg/ t^{N+1}\right)\Bigg|^2\widetilde{T}(g/ t^{N+1}).\]
Following the classical approach to shifted prime differences in $\Z$, one might take $\widetilde{T}$ to be the Fourier transform of the indicator function of the shifted irreducibles, in which case the expression above would encode information about the number of shifted irreducible differences in $A$ by the orthogonality relation.  

It turns out that after painstakingly choosing the coefficients $a_f$, one can quickly and easily gain a power saving bound on $D_q(N)$.

\begin{proof}[Proof of upper bound in Theorem \ref{MainTheorem} assuming Proposition \ref{vdc}]

Let $N\in \N$, and let $T$ be given as in Proposition \ref{vdc}, and define $\widetilde{T}$ as above. Suppose $A\subseteq\{f\in \F: \deg f \leq N\}$ with 
\[(A-A)\cap \{P-1: P\in\F\}=\emptyset.\]
First, observe that
\begin{equation}\label{eq:T(g/t^{N+1})}
\sum_{g\in \F\slash (t^{N+1})}\left|\sum_{\langle f \rangle \leq \widehat{N}}\1_A(f)e\left(fg\slash t^{N+1}\right)\right|^2T(g\slash t^{N+1})\geq |A|^2.    
\end{equation}
This follows immediately by picking up the term corresponding to $g=0$ and using the fact that $T(0)=1$. On the other hand, by applying the orthogonality relation and Plancherel, we have 
\begin{align*}
&\sum_{g\in \F\slash (t^{N+1})}\Bigg|\sum_{\langle f \rangle \leq \widehat{N}}\1_A(f)e\left(fg/ t^{N+1}\right)\Bigg|^2\widetilde{T}(g/ t^{N+1})
\\
&=\sum_{g\in \frac{\F}{\left(t^{N+1}\right)}}\sum_{\substack{\langle f_i\rangle \leq \widehat{N}\\i\in\{1,2\}}}\1_A(f_1)\1_A(f_2)e\left(\frac{(f_1-f_2)g}{t^{N+1}}\right)\Bigg(a_0+\sum_{\langle P\rangle \leq \widehat{N}}a_{P-1}e\left(\frac{g(P-1)}{t^{N+1}}\right)\Bigg)\\
&=q\widehat{N}|A|a_0+\sum_{\substack{\langle f_i\rangle \leq \widehat{N}\\i\in\{1,2\}}}\sum_{
\substack{\langle P\rangle \leq \widehat{N}}}a_{P-1}\1_A(f_1)\1_A(f_2)\sum_{g\in \F/ (t^{N+1})}e\left(\frac{(f_1-f_2+(P-1))g}{t^{N+1}}\right)\\
&=q\widehat{N}|A|a_0,
\end{align*}
where we also used the fact that $A$ has no shifted irreducible differences. Taking the real parts of both sides of this equation and comparing it with \eqref{eq:T(g/t^{N+1})}, we obtain
$|A|^2\leq q\widehat{N}|A|a_0.$
Since $0<a_0 \ll_{\epsilon,q} \widehat{N}^{-1/12+\epsilon}$, we conclude that 
$|A|\ll_{q,\epsilon} \widehat{N}^{11/12+\epsilon}$
as needed. 
\end{proof}

In practice, it helps to prove the following proposition instead, which immediately implies Proposition \ref{vdc}. This is essentially just a normalized version of $T$ from Proposition \ref{vdc}. Let $L=N(1\slash 12-\eta)$ so that $\widehat{L}=\widehat{N}^{1\slash 12-\eta}$, where $\eta$ is a small positive constant to be chosen later.

\begin{prop}\label{normalizedvdc}
Let $N\in \N$ and $\epsilon>0$. There exists a function $\Psi:\F\to \R$ which satisfies the following three properties. 

\begin{center}
\begin{tabular}{rl}
(1) & $\text{supp}\Psi\subseteq\{P-1: P\in \F\}$, \\
(2) & $\displaystyle \sum_{\langle f \rangle \leq \widehat{N}}\Psi(f)\Re[e(fx)]\geq -\frac{\widehat{N}}{\widehat{L}}$ for every $x \in K_\infty$, \\
(3) & $\displaystyle\sum_{\langle f \rangle \leq \widehat{N}}\Psi(f)\gg\widehat{N}$. \\
\end{tabular}
\end{center}
\end{prop}
Then, it is not hard to see that Proposition \ref{normalizedvdc} implies Proposition \ref{vdc} by taking 
\[T(x)=\left(\frac{\widehat{N}}{\widehat{L}}+\sum_{\langle f \rangle \leq \widehat{N}}\Psi(f) \right)^{-1}\left(\frac{\widehat{N}}{\widehat{L}}+\sum_{\langle f \rangle \leq \widehat{N}}\Psi(f)\Re[e(fx)]\right).\]

We briefly describe a function field analogue of Green's heuristic. The very rough idea is to take $\Psi(f)$ to be something like 
\[\Psi_0(f)=\Lambda'(f+1)\Lambda'(f-1)\tau(f)^2,\] 
which is clearly supported on the shifted irreducibles $P-1$ and hence satisfies (1). By the function-field analogue of the twin prime conjecture for $\F$, which was established in the remarkable work of Sawin and Shusterman \cite{SawinShusterman2022}, we have, at least when $q$ is suitably large compared to the characteristic $p>2$ of $\mathbb{F}_q$, that
\[\sum_{\langle f \rangle \leq \widehat{N}}\Psi_0(f)\ge\sum_{\langle f \rangle \leq \widehat{N}}\Lambda'(f+1)\Lambda'(f-1)\gg \widehat{N},\]
which is precisely (3). To see further why this heuristic is useful, we take $x=r\slash P$ in the statement of Proposition \ref{normalizedvdc}, where $P$ is irreducible and $r\in\F$ satisfies $(r,P)=1$. We might expect the values of $\Psi_0(f)$ to be equidistributed across the congruence classes 
\[\{g\in \F\slash(P): g\not\equiv -1,1,0\pmod P\},\] 
to have no weight at $\pm 1\pmod{P}$, and to have $4$ times as much weight at $0\pmod{P}$ as any other congruence class (because of the factor of $\tau(f)^2$). By splitting over congruence classes mod $P$, these observations might imply that the exponential sum 
\[\sum_{\langle f\rangle \leq \widehat{N}} \Psi_0(f)e(fr\slash P)=\sum_{a\in\F/(P)}e(ar/P)\sum_{\substack{\langle f\rangle\leq \widehat{N}\\f\equiv a\hspace*{-2.5mm}\pmod{P}}}\Psi_0(f)\]
is nonnegative and real-valued (by the symmetry in the definition of $\Psi_0$), which would hopefully give the property (2) in Proposition \ref{normalizedvdc}. It is natural to expect the same heuristic to hold with $P$ replaced by a product of distinct irreducibles coprime to $r$. Thus, one might hope that $\Psi_0$ satisfies Proposition \ref{normalizedvdc} at $x=\frac{r}{P}$, and perhaps this calculation would then extend to fractions with denominators which are not square free and also to all of $\mathbb{T}$. However, making this heuristic rigorous would require extremely precise information about the distribution of twin irreducibles across congruence classes in $\F$, which may be possible by leveraging ideas from \cite{SawinShusterman2022}. For a more detailed explanation of this heuristic, see either of Green's papers \cite{Green23,Green}. 

To circumvent this issue, we avoid working with $\Psi_0$, and instead construct a weight function $\Psi$ which roughly has the properties described in the heuristic above. We must take approximations of $\Lambda'$ and $\tau$ which are more amenable to analysis. We begin by defining local factors. Let $P \in \F$ be monic and irreducible. Following \cite{Green}, we define 
\[
\lambda_P(f):=
\begin{cases}
0 &\text{ if } P \mid f,\\[6pt]
\frac{\langle P \rangle }{\langle P \rangle - 1} &\text{ if } P \nmid f,
\end{cases}
\]
and 
\[
h_P(f):=
\begin{cases}
\frac{4\langle P \rangle}{\langle P \rangle + 3} &\text{ if } P \mid f,\\[6pt]
\frac{\langle P \rangle}{\langle P \rangle + 3} &\text{ if } P \nmid f.
\end{cases}
\]
Defined this way, both $\lambda_P(f)$ and $h_P(f)$ have mean value 1 as $f$ ranges over $\F/(P)$. By the orthogonality relation, we have the Fourier expansions for $\lambda_P$ and $h_P$:
\begin{align*}
\lambda_P(f) &= 1 - \frac{1}{\langle P \rangle - 1}\sum_{g \in (\F/(P))^\times} e(fg / P),\\
h_P(f) &= 1 + \frac{3}{\langle P \rangle + 3}\sum_{g \in (\F/(P))^\times} e(fg / P).
\end{align*}
Let $X>0$ and define 
\begin{align*}
\widetilde{\Lambda}_X(f)&:=\prod_{\substack{\langle P \rangle\leq \widehat{X}\\P\text{~monic, irreducible}}}\lambda_P(f),\\
\widetilde{H}_X(f)&:=\prod_{\substack{\langle P \rangle\leq \widehat{X}\\P\text{~monic, irreducible}}}h_P(f).
\end{align*}
Since $\lambda_P(f)=1_{P\nmid f}\langle P \rangle/\phi(P)$, we have
\[\widetilde{\Lambda}_X(f)=\frac{\langle \mathscr{P}_X \rangle}{\phi(\mathscr{P}_X)}1_{(f,\mathscr{P}_X)=1},\]
where $\mathscr{P}_X$ is the product of all monic irreducibles $P\in\F$ of degree at most $X$. In particular, $\widetilde{\Lambda}_X$ may be thought of as the normalized indicator function of ``$X$-rough polynomials" with mean value 1 over $\F/(\mathscr{P}_X)$, which, from a sieve-theoretic viewpoint, should approximate $\Lambda$ on $f$ of degree at most $N$ when $X$ is chosen to be a small constant multiple of $N$. Similarly, we have
\[\widetilde{H}_X(f)=\tau((f,\mathscr{P}_X))^2\prod_{\substack{\langle P \rangle\leq \widehat{X}\\P\text{~monic, irreducible}}}\left(1+\frac{3}{\langle P \rangle}\right)^{-1},\]
which may be thought of as a normalized approximation to $\tau(f)^2$.

By our Fourier expansions for $\lambda_P$ and $h_P$, we have the Ramanujan expansions of $\widetilde{\Lambda}_X(f)$ and $\widetilde{H}_X(f)$:
\begin{align*}
\widetilde{\Lambda}_X(f)&=\sum_{g\mid \mathscr{P}_X}\alpha(g)c_{g}(f)=\sum_{g\mid \mathscr{P}_X}\alpha(g)\sum_{a\in (\F\slash (g))^\times}e(af\slash g),\\
\widetilde{H}_X(f)&=\sum_{g\mid \mathscr{P}_X}\alpha'(g)c_{g}(f)=\sum_{g\mid \mathscr{P}_X}\alpha'(g)\sum_{a\in (\F\slash (g))^\times}e(af\slash g),
\end{align*}
where $\alpha$ and $\alpha'$ are supported on square-free $g$ with 
\begin{align*}
\alpha(g)&=\prod_{\substack{P\mid g\\P\text{~monic, irreducible}}}\frac{-1}{\langle P \rangle -1}=\frac{\mu(g)}{\phi(g)},\\
\alpha'(g)&=\prod_{\substack{P\mid g\\P\text{~monic, irreducible}}}\frac{3}{\langle P \rangle +3}.
\end{align*}
For technical reasons, we will work with the truncated variants of $\widetilde{\Lambda}_X$ and $\widetilde{H}_X$ given by 
\begin{align*}
\Lambda_X(f)&:=\sum_{\langle g \rangle \leq \widehat{X}}\alpha(g)\sum_{a\in (\F\slash (g))^\times}e(af\slash g),\\  
H_X(f)&:=\sum_{\langle g \rangle \leq \widehat{X}}\alpha'(g)\sum_{a\in (\F\slash (g))^\times}e(af\slash g).
\end{align*}
These truncated variants are designed to mimic $\Lambda$ and $\tau^2$, respectively. They have the advantage of being shorter sums in comparison with $\widetilde{\Lambda}_X$ and $\widetilde{H}_X$, especially considering the presence of the factor $\mu(g)$ in $\alpha(g)$. Moreover, as one might expect, they enjoy essentially the same equidistribution properties in arithmetic progressions as $\Lambda$ and $\tau^2$.

Finally, set $\widehat{Q}=\widehat{N}^{1\slash 12}$. By our discussion above, we expect
\[\Psi(f):=\Lambda'(f+1)\Lambda_Q(f-1)H_Q(f)\]
to satisfy the conclusion of Proposition \ref{normalizedvdc}. The goal of the rest of the manuscript is to show this is the case. 


\section{Outline of the proof}
Recall the definitions of $Q$ and $\Psi$ from the previous section. Also define $\widehat{R}=\widehat{N}^{1\slash 4}$. Because of the difficult nature of $\Lambda'$, we do not work directly with $\Psi$ but instead with an approximation of $\Psi$, which we call $\Psi': \F\to \R$ given by 
\[\Psi'(f):=\Lambda_R(f+1)\Lambda_Q(f-1)H_Q(f).\]
In a later section, we show $\Psi$ and $\Psi'$ are uniformly close in the Fourier sense. This way, the properties we prove for the exponential sum $\Psi'$ will readily transfer to $\Psi$. This is the content of the following proposition. (Recall that $\widehat{L}=\widehat{N}^{1\slash 12-\eta}$, where $\eta>0$ is a small constant to be chosen later.)
\begin{prop}\label{fourierclose}
Let $N\in \N$. Then, 
\[\sup_{x\in \mathbb{T}} \left|\sum_{\langle f\rangle \leq \widehat{N}} (\Psi(f)-\Psi'(f))e(fx)\right|\leq \frac{\widehat{N}}{2\widehat{L}}.\]
\end{prop}

We claim the following two estimates for all $x\in \mathbb{T}$. These correspond to properties (2) and (3) in Proposition \ref{normalizedvdc}, respectively: 

\begin{equation}\label{nottoonegative}
\sum_{\langle f\rangle \leq \widehat{N}} \Psi'(f)\Re e(fx)\geq -\frac{\widehat{N}}{2\widehat{L}}
\end{equation}
and 
\begin{equation}\label{bigatzero}
\sum_{\langle f\rangle \leq \widehat{N}} \Psi'(f)\gg \widehat{N}.
\end{equation}

With Proposition \ref{fourierclose} in hand, it suffices to prove \eqref{nottoonegative} and \eqref{bigatzero}, and then properties (2) and (3) in Proposition \ref{normalizedvdc} will follow immediately by the triangle inequality. 

In order to establish \eqref{nottoonegative} and \eqref{bigatzero}, we need to use the circle method. For this, we need to define the major and minor arcs, which determine which $x\in \mathbb{T}$ are close or far from rationals with small denominator. Let $\epsilon$ be a small positive constant to be chosen later. Let $\mathfrak{M}$ denote the set of $x\in \mathbb{T}$ such that there exist coprime $r,b \in \F$ with 
\begin{equation}\label{majorarcparams}
x= -\frac{b}{r}+\eta,\quad \langle r \rangle \ll_{\epsilon} \widehat{L}^{1+\epsilon},\quad\langle \eta \rangle \ll_{\epsilon} \widehat{L}^{1+\epsilon}\widehat{N}^{-1}\langle r \rangle^{-1}.
\end{equation}
 Let $\mathfrak{m}:=\mathbb{T}\setminus\mathfrak{M}$. We say that $\mathfrak{M}$ is the set of \textit{major arcs} and $\mathfrak{m}$ is the set of \textit{minor arcs}. 
\begin{prop}\label{minorarcs} Let $N\in \N$ and $x\in \mathbb{T}$. If 
\[\left|\sum_{\langle f\rangle \leq \widehat{N}} \Psi'(f)\Re e(fx)\right|>\frac{\widehat{N}}{2\widehat{L}},\]
then $x\in \mathfrak{M}$. 
\end{prop}
The point of Proposition \ref{minorarcs} is that if \eqref{nottoonegative} fails, then $x$ lies in a major arc, in which case we have a way of nicely approximating the exponential sum. Indeed, consider the coefficients $\beta^{\text{trunc}}(\lambda)$ for $\lambda\in K$ given by 
\[\Lambda_R(f+1)\Lambda_Q(f-1)H_Q(f)=\sum_{\lambda \in K\slash \F} \beta^{\text{trunc}}(\lambda) e(\lambda f). \]
Then, we have the following formula for the exponential sum of $\Psi'$ when $x\in \mathfrak{M}$.

\begin{prop}\label{majorarcformula}
Let $N\in \N$. If $x\in \mathfrak{M}$ and $b,r,\eta$ are given as in \eqref{majorarcparams}, then 
\[\sum_{\langle f\rangle \leq \widehat{N}} \Psi'(f)e(fx) =\beta^{\text{trunc}}(b\slash r)\sum_{\langle f\rangle \leq \widehat{N}} e(f\eta).\]
\end{prop}

Unfortunately, $\beta^{\text{trunc}}$ is still hard to work with, so we need another layer of approximation. Consider the coefficients $\beta(\lambda)$ for $\lambda\in K$ given by 
\[\widetilde{\Lambda}_R(f+1)\widetilde{\Lambda}_R(f-1)\widetilde{H}_R(f)=\sum_{\lambda \in K\slash \F} \beta(\lambda) e(\lambda f). \]
We shall show that $\beta^{\text{trunc}}(\lambda)$ and $\beta(\lambda)$ are uniformly close for $\lambda$ with small denominator. This in turn will allow us to just consider $\beta$, which turns out to be much nicer. 
\begin{prop}\label{betasclose}Suppose $\lambda\in K$ with $\langle \text{denom}(\lambda)\rangle \ll \widehat{L}^{1+\epsilon}$. Then 
\[|\beta^{\text{trunc}}(\lambda)-\beta(\lambda)|\leq \frac{1}{10\widehat{L}}.\]
\end{prop}

Propositions \ref{majorarcformula} and \ref{betasclose} combined yield instantly the following estimate for the exponential sum of $\Psi'$, with which we will be able to quickly derive \eqref{nottoonegative} and \eqref{bigatzero}. 

\begin{prop}\label{majorarcestimate}
Let $N\in \N$. If $x\in \mathfrak{M}$ with $b,r,\eta$ given as in \eqref{majorarcparams}, then
\[\left|\sum_{\langle f\rangle \leq \widehat{N}} \Psi'(f)\Re e(fx)-\beta(b\slash r) \sum_{\langle f \rangle \leq \widehat{N}}\Re e(f\eta) \right|\leq \frac{\widehat{N}}{2\widehat{L}}.\]
\end{prop}

Finally, we establish two crucial properties of $\beta$: it is real and nonnegative and $\beta(0)\gg 1$. 

\begin{prop}\label{realandnonnegative}
For every $\lambda \in K\slash \F$, $\beta(\lambda)$ is real and nonnegative, and $\beta(0)\gg 1$. 
\end{prop}
\begin{proof}
For each monic, irreducible $P\in \F$, define $u_P:\F\to \R$ by 
\begin{align*}
u_P(f)&=\left(1-\frac{1}{ \langle P \rangle}\right)^2\left(1+\frac{3}{\langle P \rangle}\right)\lambda_P(f-1)\lambda_P(f+1)h_P(f)\\
&=
\begin{cases}
0 &\text{ if } f\equiv \pm 1\pmod{P},\\
4 &\text{ if } f\equiv 0\pmod{P},\\
1 &~\text{otherwise}.
\end{cases}
\end{align*}
Keeping in mind that $\text{char}(\mathbb{F}_q)>2$, we have
\begin{align*}
\widehat{u}_P(g)&=\frac{1}{\langle P \rangle}\sum_{f\in \F\slash (P)}u_P(f)e\left(\frac{-fg}{P}\right)\\
&=
\begin{cases}
\frac{1}{\langle P \rangle}\big(3-2\Re e(g\slash P)\big) &\text{ if } g\not\equiv 0\pmod{P},\\
1+\frac{1}{\langle P \rangle}&\text{ if } g\equiv 0\pmod{P},
\end{cases}
\end{align*}
so, by Fourier inversion, 
\[u_P(f)=1+\frac{1}{\langle P \rangle}+\sum_{g\in (\F\slash (P))^\times}\frac{1}{\langle P \rangle}\left(3-2\Re e(g\slash P)\right)e(fg\slash P).\]
Now, 
\[\tilde{\Lambda}_R(f+1)\tilde{\Lambda}_R(f-1)\tilde{H}_R(f)=\prod_{\substack{\langle P \rangle \leq \widehat{R}\\ P~\text{monic, irreducible}}}\left(1-\frac{1}{ \langle P \rangle}\right)^{-2}\left(1+\frac{3}{\langle P \rangle}\right)^{-1}u_P(f).\]
Substituting our Fourier expansions for $u_P$ into this formula and combining terms, we see that $\beta(\lambda)$ is real and nonnegative, since the Fourier coefficients in the expansion of each $u_P$ are real and nonnegative. Moreover, upon inspection, we have 
\[\beta(0)\geq \prod_{\substack{\langle P \rangle \leq \widehat{R}\\ P~\text{monic, irreducible}}}\left(1-\frac{1}{ \langle P \rangle}\right)^{-2}\left(1+\frac{3}{\langle P \rangle}\right)^{-1}\left(1+\frac{1}{\langle P \rangle }\right)\gg1,\]
as needed.
\end{proof}

With all of these propositions in hand, we may prove \eqref{nottoonegative} and \eqref{bigatzero}, which, when combined with Proposition \ref{fourierclose}, imply Proposition \ref{normalizedvdc}, as noted above.

\begin{proof}[Proofs of \eqref{nottoonegative} and \eqref{bigatzero} assuming the propositions above]

By Proposition \ref{minorarcs}, if \eqref{nottoonegative} fails, then $x$ satisfies \eqref{majorarcparams}. Thus, by Proposition \ref{majorarcestimate}, we have 
\[\sum_{\langle f\rangle \leq \widehat{N}} \Psi'(f)\Re e(fx) \geq \beta(b\slash r) \sum_{\langle f \rangle \leq \widehat{N}}\Re e(f\eta)-\frac{\widehat{N}}{2\widehat{L}}\geq -\frac{\widehat{N}}{2\widehat{L}}.\]
Here, we used the fact that 
\[\beta(b\slash r) \sum_{\langle f \rangle \leq \widehat{N}}\Re e(f\eta)\]
is real and nonnegative by Proposition \ref{realandnonnegative} and Lemma \ref{exponentialsumformula}. This confirms \eqref{nottoonegative}. 

To prove \eqref{bigatzero}, notice first that $0\in\mathfrak{M}$ trivially. Thus, by Proposition \ref{majorarcestimate}, Proposition \ref{realandnonnegative} and Lemma \ref{exponentialsumformula}, we have  
\[\sum_{\langle f\rangle \leq \widehat{N}} \Psi'(f) \geq \beta(0) q\widehat{N}-\frac{\widehat{N}}{2\widehat{L}}\gg \widehat{N},\]
as needed. 
\end{proof}

The remaining parts of the paper will be dedicated to the proofs of Propositions \ref{fourierclose}--\ref{betasclose}.

\section{Products of Ramanujan expansions}
Before proving Propositions \ref{fourierclose}--\ref{betasclose}, we need some preparatory results concerning products of Ramanujan expansions. Define $\mathcal{C}_B(X)$ to be the set of exponential sums 
\[S(f)=\sum_{\lambda\in K\slash \F}c(\lambda)e(f\lambda)\]
where $c(\lambda)$ is supported on those $\lambda= a\slash s$ where $s$ is square-free with $\langle s \rangle \leq \widehat{X}$, and $|c(\lambda)|\leq \tau^B(s)\slash \langle s\rangle $. These can be compared to classical Ramanujan--Fourier expansions \cite{Ramanujan1918}. We will show that $\Lambda_X, H_X\in \mathcal{C}_2(X)$, and we will argue that functions from classes $\mathcal{C}_B(X)$ are weakly closed under multiplication, so that, in particular, $\Psi'(f)=\Lambda_R(f+1)\Lambda_Q(f-1)H_Q(f)$ will also be of this form. 

We need the following elementary estimate. 

\begin{lem}\label{denomcount}
For any $d,r,s\in\F\setminus\{0\}$ and $b\in\F$ with $(b,r)=1$, we have
\[|\{a\in \F\slash (s): \text{denom}(a\slash s+b\slash r)=d\}|\leq \frac{\langle d \rangle \langle (r,s)\rangle }{\langle r \rangle}.\]
\end{lem}
\begin{proof}
Write $r'=r/(r,s)$ and $s'=s/(r,s)$. Since $[r,s]=rs/(r,s)$, we have
\[\frac{a}{s}+\frac{b}{r}=\frac{ar+bs}{rs}=\frac{ar'+bs'}{[r,s]}.\]
If $\text{denom}(a\slash s+b\slash r)=d$, then $d\mid [r,s]$, and $ar'+bs'$ must be a multiple of $h:=[r,s]/d$. In other words, we have $ar'+bs'\equiv 0\pmod{h}$. It is easy to see that $(r',h)=1$. Indeed, if $g\mid(r',h)$ for some irreducible $g\in\F$, then $g\mid bs'$, which would imply $g\mid s'$, due to our assumption that $(b,r)=1$. But this is impossible, since $(r',s')=1$. As a consequence, there exists a unique $x\in \F\slash (h)$ satisfying $xr'+bs'\equiv 0\pmod{h}$. Since $x$ has exactly $\langle s \rangle/\langle h \rangle$ lifts in $\F\slash (s)$, we conclude 
\[|\{a\in \F\slash (s): \text{denom}(a\slash s+b\slash r)=d\}|\le \frac{\langle s \rangle}{\langle h \rangle} =\frac{\langle d \rangle \langle (r,s)\rangle }{\langle r \rangle},\]
as required.
\end{proof}

The following lemma says that $\mathcal{C}_B(X)$ is weakly closed under multiplication, as alluded to above. 
\begin{lem}\label{exposumclosed}
If $S_1\in \mathcal{C}_{B_1}(X_1)$ and $S_2\in \mathcal{C}_{B_2}(X_2)$, then 
\[S_1S_2\in \big(\log^{O_{B_1,B_2}(1)}\widehat{X_1}\big)\mathcal{C}_{B_1+2B_2+3}(X_1+X_2).\]
\end{lem}
\begin{proof}
Suppose 
\begin{align*}
S_1(f)&=\sum_{\lambda\in K\slash \F}c_1(\lambda)e(f\lambda),\\
S_2(f)&=\sum_{\lambda\in K\slash \F}c_2(\lambda)e(f\lambda),\\
S_1S_2(f)&=\sum_{\lambda\in K\slash \F}c(\lambda)e(f\lambda).
\end{align*}
Then,
\begin{align*}
|c(a\slash s)|&=\Bigg|\sum_{\substack{\langle s_1\rangle \leq \widehat{X_1}\\\langle s_2\rangle \leq \widehat{X_2}}}\sum_{\substack{a_1\in (\F\slash (s_1))^\times \\a_2\in (\F\slash (s_2))^\times\\\frac{a_1}{s_1}+\frac{a_2}{s_2}=\frac{a}{s}}} c_1(a_1\slash s_1)c_2(a_2\slash s_2)\Bigg|\\
&\leq \sum_{\substack{\langle s_1 \rangle \leq \widehat{X_1}\\\mu^2(s_1)=1}}\frac{\tau(s_1)^{B_1}}{\langle s_1 \rangle }\sum_{d\mid s_1s}\frac{\tau(d)^{B_2}}{\langle d \rangle }\bigg|\bigg\{a_1\in \bigg(\frac{\F}{(s_1)}\bigg)^\times: \text{denom}\bigg(\frac{a_1}{s_1}-\frac{a}{s}\bigg)=d\bigg\}\bigg|.
\end{align*}
Applying Lemma \ref{denomcount}, this is bounded above by 
\[\sum_{\substack{\langle s_1 \rangle \leq \widehat{X_1}\\\mu^2(s_1)=1}}\frac{\tau(s_1)^{B_1}}{\langle s_1 \rangle }\sum_{d\mid s_1s}\tau(d)^{B_2}\frac{\langle (s_1,s)\rangle }{\langle s \rangle }\leq \frac{\tau(s)^{B_2+1}}{\langle s \rangle }\sum_{\substack{\langle s_1 \rangle \leq \widehat{X_1}\\\mu^2(s_1)=1}}\frac{\tau(s_1)^{B_1+B_2+1}\langle(s_1,s)\rangle }{\langle s_1\rangle },\]
where we used the trivial bound 
\[\sum_{d\mid s_1s}\tau(d)^{B_2}\leq \tau(s_1)^{B_2+1}\tau(s)^{B_2+1}.\]
Then the lemma follows by Lemma \ref{rankinbound} in the appendix.
\end{proof}

Finally, we apply Lemma \ref{exposumclosed} to the Ramanujan series $\Lambda_X$ and $H_X$. Recall that 
\begin{align*}
\Lambda_X(f)&=\sum_{\langle g \rangle \leq \widehat{X}}\alpha(g)\sum_{a\in (\F\slash (g))^\times}e(af\slash g),\\
H_X(f)&=\sum_{\langle g \rangle \leq \widehat{X}}\alpha'(g)\sum_{a\in (\F\slash (g))^\times}e(af\slash g),
\end{align*}
where 
\[\alpha(g)=\frac{\mu(g)}{\phi(g)}\]
and 
\[\alpha'(g)=\prod_{\substack{P\mid g\\P\text{~monic, irreducible}}}\frac{3}{\langle P \rangle +3}.\]
It is not hard to see that 
\[|\alpha(g)|,|\alpha'(g)|\leq \frac{\tau(g)^2}{\langle g \rangle}.\]
Thus, we have $\Lambda_R\in \mathcal{C}_2(R)$ and $\Lambda_Q, H_Q\in \mathcal{C}_2(Q)$, so 
\[\Lambda_Q(f-1)H_Q(f)\in \big(\log^{O(1)}\widehat{N}\big)\mathcal{C}_9(N\slash 6)\]
and 
\[\Psi'(f)=\Lambda_R(f+1)\Lambda_Q(f-1)H_Q(f)\in \big(\log^{O(1)}\widehat{N}\big)\mathcal{C}_{16}(5N\slash 12).\]

\section{Passing to a proxy for the prime polynomials}\label{sec:proof of Prop 5.1}
The goal of this section is to prove Proposition \ref{fourierclose}, which says that $\Psi$ and $\Psi'$ are uniformly close in the Fourier sense. The key to accomplishing this goal is the following result which shows that $\Lambda'$ and $\Lambda_R$ are uniformly close in the Fourier sense. 

\begin{prop}\label{Lambdafourierclose} For any $N\in \N$ and $\epsilon>0$, we have
\[\sup_{x\in \mathbb{T}} \Bigg|\sum_{\langle f\rangle \leq \widehat{N}} (\Lambda'(f)-\Lambda_R(f))e(fx)\Bigg|\ll_\epsilon \widehat{N}^{3\slash 4+\epsilon}.\]
\end{prop}

Assuming Proposition \ref{Lambdafourierclose}, we are prepared to prove Proposition \ref{fourierclose}.

\begin{proof}[Proof of Proposition \ref{fourierclose} assuming Proposition \ref{Lambdafourierclose}]

First, we apply Lemma \ref{exposumclosed} to write 
\[\Lambda_Q(f-1)H_Q(f)=\sum_{\substack{\lambda\in K\slash \F\\ \langle \text{denom} (\lambda)\rangle\leq \widehat{N}^{1\slash 6}}} c(\lambda)e(\lambda f),\]
with $|c(\lambda)|\leq \frac{\tau(s)^9}{\langle s\rangle }$, where $s=\text{denom}(\lambda)$. 
By Proposition \ref{Lambdafourierclose} and the definitions of $\Psi$ and $\Psi'$, we have 
\begin{align*}
&\Bigg|\sum_{\langle f\rangle \leq \widehat{N}} (\Psi(f)-\Psi'(f))e(fx)\Bigg|\\
&=\Bigg|\sum_{\substack{\lambda\in K\slash \F\\ \langle \text{denom} (\lambda)\rangle\leq \widehat{N}^{1\slash 6}}} c(\lambda)\sum_{\langle f\rangle \leq \widehat{N}} (\Lambda'(f+1)-\Lambda_R(f+1))e(f(x+\lambda))\Bigg|\\
&\ll_\epsilon \widehat{N}^{3\slash 4+\epsilon}\sum_{\substack{\lambda\in K\slash \F\\ \langle \text{denom} (\lambda)\rangle\leq \widehat{N}^{1\slash 6}}} |c(\lambda)|.
\end{align*}
Finally, by the triangle inequality and trivial standard divisor bound $\tau(s)\ll_\epsilon \langle s \rangle ^\epsilon$, we have
\begin{align*}
\sum_{\substack{\lambda\in K\slash \F\\ \langle \text{denom} (\lambda)\rangle\leq \widehat{N}^{1\slash 6}}} |c(\lambda)|&\leq \sum_{\substack{\lambda\in K\slash \F\\ \langle \text{denom} (\lambda)\rangle\leq \widehat{N}^{1\slash 6}}} \frac{\tau(s)^9}{\langle s\rangle }\\
&\ll_\epsilon \widehat{N}^{\epsilon}\sum_{\substack{\lambda\in K\slash \F\\ \langle \text{denom} (\lambda)\rangle\leq \widehat{N}^{1\slash 6}}} \frac{1}{\langle s\rangle }\\
&\ll_\epsilon \widehat{N}^{1\slash 6+\epsilon}.
\end{align*}
Combining these estimates gives 
\[\Bigg|\sum_{\langle f\rangle \leq \widehat{N}} (\Psi(f)-\Psi'(f))e(fx)\Bigg|\leq C_\epsilon \widehat{N}^{3\slash 4+\epsilon}\widehat{N}^{1\slash 6+\epsilon}\leq \frac{\widehat{N}}{2\widehat{L}},\]
upon taking $\epsilon$ sufficiently small and $N$ sufficiently large. 
\end{proof}

We now turn to the proof of Proposition \ref{Lambdafourierclose}. For any $x\in \mathbb{T}$, Dirichlet's approximation theorem yields coprime $b,r \in \F$ with 
\begin{equation}\label{dirichletx}
x=b\slash r+\eta, \quad\langle \eta \rangle <(\langle r \rangle \widehat{N}^{1\slash 2})^{-1}, \quad\langle r \rangle \leq \widehat{N}^{1\slash 2}.
\end{equation}
We need the following exponential sum estimate for the von Mangoldt function from a paper of Hayes on the ternary Goldbach problem in function fields \cite{Hayes66}. Due to Weil's proof of GRH for $\F$, this estimate is now known to hold unconditionally. 


\begin{thmF}[Hayes, Weil]\label{Hayes}
Let $n\leq N$ be a positive integer and suppose $x=b\slash r +\eta$ with $b,r,\eta$ satisfying \eqref{dirichletx}. Then for any $\epsilon>0$, 
\[\Bigg|\sum_{\substack{f~\text{monic}\\\deg f=n}}\Lambda'(f)e(fx)\;-\;\frac{\mu(r)}{\phi(r)}\sum_{\substack{f~\text{monic}\\ \deg f=n}}e(f\eta)\Bigg|\ll_\epsilon \widehat{N}^{3\slash 4+\epsilon}.\]
\end{thmF}
\begin{proof}
If $\langle r \rangle \leq \widehat{n}^{1/2}$, then the result follows immediately from \cite[Lemma 5]{Pollack13} and Corollary \ref{exponentialsumformulaV2}. 

Suppose now that $\langle r \rangle >\widehat{n}^{1/2}$. Since \eqref{dirichletx} implies that $n<N$ and
\[\text{ord}(\eta)<-\deg r-\frac{N}{2}\le -\left(\left\lfloor \frac{n}{2}\right\rfloor+1\right)-\frac{n+1}{2}\le -n-1,\] 
we have $e(f\eta)=1$ for all monic $f\in \F$ with $\deg f=n$. Thus, 
\[\sum_{\substack{f~\text{monic}\\\deg f=n}}\Lambda'(f)e(fx)=\sum_{\substack{f~\text{monic}\\\deg f=n}}\Lambda'(f)e\left(\frac{fb}{r}\right).\]
By the orthogonality relations for Dirichlet characters, this becomes 
\begin{align*}
&\sum_{\substack{f~\text{monic}\\\deg f=n\\(f,r)=1}}\Lambda'(f)e\left(\frac{fb}{r}\right)+O(n\deg r)\\
&=\sum_{g\in (\F\slash (r))^\times}e\left(\frac{gb}{r}\right)\sum_{\substack{f~\text{monic}\\\deg f=n\\f\equiv g\hspace*{-2.5mm}\pmod{r}}}\Lambda'(f)+O(n\deg r)\\
&=\frac{1}{\phi(r)}\sum_{g\in (\F\slash (r))^\times}e\left(\frac{gb}{r}\right)\sum_{\substack{f~\text{monic}\\\deg f=n}}\Lambda'(f)\sum_{\chi\hspace*{-2.5mm}\pmod{r}}\bar{\chi}(g)\chi(f)+O(n\deg r)\\
&=\frac{1}{\phi(r)}\sum_{\chi\hspace*{-2.5mm}\pmod{r}}\sum_{g\in (\F\slash (r))^\times}\bar{\chi}(g)e\left(\frac{gb}{r}\right)\sum_{\substack{f~\text{monic}\\\deg f=n}}\Lambda'(f)\chi(f)+O(n\deg r)\\
&=\frac{1}{\phi(r)}\sum_{\chi\hspace*{-2.5mm}\pmod{r}}\chi(b)G(\bar{\chi})\psi(n,\chi)+O(n\deg r).
\end{align*}
where $\psi(n,\chi)$ is defined as in \eqref{eq:psi(n,chi)}. The contribution to the last line above from the trivial character $\chi_0$ is
\[\frac{\mu(r)}{\phi(r)}\sum_{\substack{f~\text{monic}\\\deg f=n\\(f,r)=1}}\Lambda'(f)=\frac{\mu(r)}{\phi(r)}q^n+O\left(\frac{q^{n/2}+n\deg r}{\phi(r)}\right)\]
by \eqref{eq:PNT}. On the other hand, the character sum estimate \eqref{eq:charactersum} in conjunction with the upper bound for the Gauss sums (Lemma \ref{lem:Gausssum}) implies that the contribution from the nontrivial characters is $O(\langle r \rangle ^{1/2}q^{n/2}(\deg r+n))$. Combining these contributions, we obtain
\[\sum_{\substack{f~\text{monic}\\\deg f=n}}\Lambda'(f)e(fx)=\frac{\mu(r)}{\phi(r)}q^n+O\left(\langle r \rangle ^{1\slash 2}\widehat{n}^{1\slash 2}(\deg r+n)\right).\]
The result follows now from Corollary \ref{exponentialsumformulaV2}, as needed. 
\end{proof}

Theorem F yields the following corollary. 

\begin{cor}\label{prop:FourierLambda'}
Let $N\in \N$ and suppose $x\in \mathbb{T}$ with $b,r,\eta$ given as in \eqref{dirichletx}. Then, for any $\epsilon>0$, we have
\[
\quad\left|
   \sum_{\langle f\rangle\leq \widehat{N} } \Lambda'(f)\,e(fx)
   \;-\;
   \frac{\mu(r)}{\phi(r)}\sum_{\langle f\rangle\le \widehat{N}}e(f\eta)
\right|
\;
\ll_\epsilon
\widehat{N}^{3\slash 4+\epsilon}.
\]
\end{cor}

Next, we show that the exponential sum of $\Lambda_R$ enjoys the same estimate, which would then imply Proposition \ref{Lambdafourierclose} by the triangle inequality. 
\begin{prop}\label{prop:FourierLambda_R}
Let $N\in \N$ and suppose $x\in \mathbb{T}$ with $b,r,\eta$ given as in \eqref{dirichletx}. Then, 
\begin{equation}
\left|
   \sum_{\langle f\rangle\leq \widehat{N} } \Lambda_R(f)\,e(fx)
   \;-\;
   \frac{\mu(r)}{\phi(r)}
   \sum_{\langle f\rangle\leq \widehat{N}} e(f\eta)
\right|
\;\ll_\epsilon\;
\widehat{N}^{3\slash 4+\epsilon}.
\end{equation}
\end{prop}

\begin{proof}

We have 
\begin{align*}
\sum_{\langle f\rangle\leq \widehat{N}} \Lambda_{R}(f)e(fx)&= \mathbf{1}_{\langle r\rangle  \leq \widehat{R}}\frac{\mu(r)}{\phi(r)} \sum_{\langle f\rangle\leq \widehat{N}} e(f\eta)\\
&\;\;\;\;\;+ \sum_{\substack{\langle h \rangle  \leq \widehat{R}, h\neq 0}}\sum_{\substack{a \in (\F\slash (h))^\times\\ \frac{a}{h} \neq -\frac{b}{r}}} \frac{\mu(h)}{\phi(h)} \sum_{\langle f\rangle\leq \widehat{N}} e\left( \left( \frac{a}{h} + \frac{b}{r} + \eta \right)f \right).
\end{align*}
However, notice that when $\langle h\rangle \leq \widehat{R}$ we have 
\[\left\langle \frac{a}{h}+\frac{b}{r}\right\rangle \geq \left\langle \frac{1}{hr}\right\rangle \geq (\langle h\rangle \langle r\rangle)^{-1}\geq (\langle r\rangle\widehat{R})^{-1}.\]
On the other hand, we have 
\[\langle \eta \rangle < (\langle r\rangle\widehat{N}^{1\slash 2})^{-1}<(\langle r\rangle\widehat{R})^{-1}.\]
Thus, 
\[\left\langle \frac{a}{h} + \frac{b}{r} + \eta \right\rangle \geq (\langle r\rangle\widehat{R})^{-1}\geq \widehat{N}^{-3\slash 4}.\]
Therefore, by Lemma \ref{exponentialsumformula}, we see 
\[ \sum_{\langle f\rangle\leq \widehat{N}} e\left( \left( \frac{a}{h} + \frac{b}{r} + \eta \right)f\right)=0,\]
so
\[\sum_{\langle f\rangle\leq \widehat{N}} \Lambda_{R}(f)e(fx)= \mathbf{1}_{\langle r\rangle  \leq \widehat{R}}\frac{\mu(r)}{\phi(r)} \sum_{\langle f\rangle\leq \widehat{N}} e(f\eta).\]
If $\langle r \rangle \leq \widehat{R}$, then the proposition is immediate. Otherwise, we have $\langle r \rangle >\widehat{R}$, in which case 
\[\sum_{\langle f\rangle\leq \widehat{N}} \Lambda_{R}(f)e(fx)=0.\]
Hence, the proposition follows from the crude bound 
\[\Bigg|\frac{\mu(r)}{\phi(r)} \sum_{\langle f\rangle\leq \widehat{N}} e(f\eta)\Bigg|\ll_\epsilon \frac{\widehat{N}}{\widehat{R}^{1-\epsilon}}=\widehat{N}^{3\slash 4+\epsilon}.\]
This proves the proposition.
\end{proof}

Now, Proposition \ref{Lambdafourierclose} follows immediately upon combining Corollary \ref{prop:FourierLambda'} with Proposition \ref{prop:FourierLambda_R} and using the triangle inequality.

\section{The major and minor arcs}

Recall our setup \eqref{majorarcparams} says that the major arcs $\mathfrak{M}$ consist of $x\in \mathbb{T}$ such that there exist coprime $r,b \in \F$ with 
\[x= -\frac{b}{r}+\eta,\quad \langle r \rangle \ll_{\epsilon} \widehat{L}^{1+\epsilon},\quad\langle \eta \rangle \ll_{\epsilon} \widehat{L}^{1+\epsilon}\widehat{N}^{-1}\langle r \rangle^{-1}.\] 
The objective of this section is to prove Propositions \ref{minorarcs} and \ref{majorarcformula}. We restate Proposition \ref{minorarcs} below, which asserts that if the exponential sum of $\Psi'$ at $x$ is large, then $x$ is in the major arcs. 

\newcounter{saveprop}
\setcounter{saveprop}{\value{prop}}

\begingroup
  \renewcommand{\theprop}{5.\arabic{prop}}

  \setcounter{prop}{1} 
  \begin{prop}\label{prop:5.2}
 Let $x\in \mathbb{T}$. If 
\[\left|\sum_{\langle f\rangle \leq \widehat{N}} \Psi'(f)\Re e(fx)\right|>\frac{\widehat{N}}{2\widehat{L}},\]
then $x\in \mathfrak{M}$. 
  \end{prop}

\begin{proof}
Suppose $x\in \mathbb{T}$ and 
\[\left|\sum_{\langle f\rangle \leq \widehat{N}} \Psi'(f)\Re e(fx)\right|>\frac{\widehat{N}}{2\widehat{L}}.\]
Then, we immediately have 
\[\left|\sum_{\langle f\rangle \leq \widehat{N}} \Psi'(f) e(fx)\right|>\frac{\widehat{N}}{2\widehat{L}}\]
as well. By Dirichlet's approximation theorem, there exist coprime $b, r\in \F$ with 
\[x=-b\slash r+\eta, \quad\langle r \rangle \leq \widehat{N}^{5\slash 12}, \quad\langle \eta \rangle < \langle r\rangle^{-1}\widehat{N}^{-5\slash 12}.\]
Recall also that 
\[\Psi'(f)=\sum_{\lambda\in K\slash \F}\beta^{\text{trunc}}(\lambda)e(\lambda f),\]
and, by Lemma \ref{exposumclosed}, 
if $\lambda = a\slash s$, then $|\beta^{\text{trunc}}(\lambda)|\leq (\log \widehat{N})^{O(1)}\frac{\tau(s)^{16}}{\langle s \rangle }$ and $\langle s \rangle \leq \widehat{N}^{5\slash 12}$. By the triangle inequality, 
\[\sum_{\langle s \rangle \leq \widehat{N}^{5\slash 12}}\sum_{a\in (\F\slash (s))^\times}(\log \widehat{N})^{O(1)}\frac{\tau(s)^{16}}{\langle s \rangle }\left|\sum_{\langle f \rangle \leq \widehat{N}}e((a\slash s -b\slash r+\eta)f)\right|>\frac{\widehat{N}}{2\widehat{L}}.\]
Suppose for the moment that $\frac{a}{s}\neq \frac{b}{r}$ and $\langle s\rangle, \langle r \rangle \leq \widehat{N}^{5\slash 12}$. Then, we have 
\[\left\langle \frac{a}{s}-\frac{b}{r}\right\rangle \geq \left\langle \frac{1}{sr}\right\rangle \geq (\langle s\rangle \langle r\rangle)^{-1}\geq \langle r \rangle^{-1}\widehat{N}^{-5\slash 12}.\]
On the other hand, we have $\langle \eta \rangle < \langle r\rangle^{-1}\widehat{N}^{-5\slash 12}$. Thus, 
\[\left\langle \frac{a}{s} - \frac{b}{r} + \eta \right\rangle \geq \widehat{N}^{-5\slash 6}.\]
Therefore, by Lemma \ref{exponentialsumformula}, we see 
\[ \sum_{\langle f\rangle\leq \widehat{N}} e\left( \left( \frac{a}{s} - \frac{b}{r} + \eta \right)f\right)=0.\]
This forces 
\[(\log \widehat{N})^{O(1)}\frac{\tau(r)^{16}}{\langle r \rangle }\left|\sum_{\langle f \rangle \leq \widehat{N}}e(\eta f)\right|>\frac{\widehat{N}}{2\widehat{L}}.\]
Thus, by Lemma \ref{exponentialsumformula} again, we must have $\text{ord}(\eta)<-N-1$, in which case the inequality becomes 
\[(\log \widehat{N})^{O(1)}\frac{\tau(r)^{16}}{\langle r \rangle }q\widehat{N}>\frac{\widehat{N}}{2\widehat{L}}.\]
Rearranging this yields 
\[\langle r \rangle \leq (\log \widehat{N})^{O(1)}\tau(r)^{16}\widehat{L}\ll_\epsilon \widehat{L}^{1+\epsilon},\]
thus completing the proof. 
\end{proof}

On the other hand, if $x$ is in the major arcs, then we have Proposition \ref{majorarcformula}, which yields the formula for the exponential sum of $\Psi'$ at $x$. Surprisingly, the proof of this is almost the same as the proof of Proposition \ref{minorarcs}. 

  \setcounter{prop}{2} 
  \begin{prop}\label{prop:5.3}
If $x\in \mathfrak{M}$ as in \eqref{majorarcparams}, then 
\[\sum_{\langle f\rangle \leq \widehat{N}} \Psi'(f)e(fx) =\beta^{\text{trunc}}(b\slash r)\sum_{\langle f\rangle \leq \widehat{N}} e(f\eta).\]
  \end{prop}

\begin{proof}
Suppose $x\in \mathfrak{M}$ such that 
\[x= -\frac{b}{r}+\eta, \quad\langle r \rangle \ll_{\epsilon} \widehat{L}^{1+\epsilon},\quad\langle \eta \rangle \ll_{\epsilon} \widehat{L}^{1+\epsilon}\widehat{N}^{-1}\langle r \rangle^{-1}.\]
By the definition of $\Psi'$ and Lemma \ref{exposumclosed}, we have
\[
\sum_{\langle f\rangle \leq \widehat{N}} \Psi'(f)e(fx)=\sum_{\substack{\lambda\in K\slash \F\\\langle \text{denom}(\lambda)\rangle\leq\widehat{N}^{5\slash 12}}}\beta^{\text{trunc}}(\lambda)\sum_{\langle f \rangle \leq \widehat{N}}e((\lambda-b\slash r+\eta)f).
\]
On the other hand, if $\lambda\neq b/r$, then Lemma \ref{exponentialsumformula} gives 
\[\sum_{\langle f \rangle \leq \widehat{N}}e((\lambda-b\slash r+\eta)f)=0.\]
Thus, 
\[\sum_{\langle f\rangle \leq \widehat{N}} \Psi'(f)e(fx)=\beta^{\text{trunc}}\left(\frac{b}{r}\right)\sum_{\langle f \rangle \leq \widehat{N}}e(\eta f),\]
as needed. 
\end{proof}
\endgroup

\setcounter{prop}{\value{saveprop}}

\section{Triple correlation of Ramanujan expansions}
Finally, we turn to the proof of Proposition \ref{betasclose}, which claims that for every $\lambda\in K$ with $\langle \text{denom}(\lambda)\rangle \ll \widehat{L}^{1+\epsilon}$, we have
\[|\beta^{\text{trunc}}(\lambda)-\beta(\lambda)|\leq \frac{1}{10\widehat{L}}.\]
The proof turns out to be quite difficult. Unfortunately, it does not seem possible to simplify the arguments in this section using the special properties of $\F$, so our argument follows Green's \cite{Green} almost verbatim. 

We will prove the following proposition, which implies Proposition \ref{betasclose} as a consequence. 

\begin{prop}\label{triplecorrelation}
Let $B,\epsilon>0$ and consider the three series 
\[S_i(f)=\sum_{h\in \F\setminus\{0\}}\alpha_i(h)\sum_{a\in (\F/(h))^\times} e\left(\frac{af}{h}\right)\]
with $|\alpha_i(h)|\leq \tau(h)^B/\langle h \rangle$ and $\alpha_i$ supported on square-free $h$ for $i\in\{1,2,3\}$. Let $y_1,y_2,y_3\in \F$ be distinct. Consider the triple correlation 
\[S_1(f+y_1)S_2(f+y_2)S_3(f+y_3)=\sum_{\lambda\in K/\F}\beta(\lambda)e(\lambda f).\]
Let $X_1,X_2,X_3>0$ and $X=\min\{X_1,X_2,X_3\}$. For each $i\in\{1,2,3\}$, write
\[S_{i,X_i}=\sum_{0<\langle h \rangle \leq \widehat{X_i}}\alpha_i(h)\sum_{a\in (\F/(h))^\times}e\left(\frac{af}{h}\right).\]
Then for the truncated correlation
\[S_{1,X_1}(f+y_1)S_{2,X_2}(f+y_2)S_{3,X_3}(f+y_3)=\sum_{\lambda\in K/\F}\beta^{\text{trunc}}(\lambda)e(\lambda f),\]
we have
\[|\beta(\lambda)-\beta^{\text{trunc}}(\lambda)|\ll_{y_1,y_2,y_3,B,\epsilon} \langle \text{denom}(\lambda)\rangle ^{\epsilon}\widehat{X}^{\epsilon-1}.\]
\end{prop}

\hspace{2cm}

To see why Proposition \ref{betasclose} follows from Proposition \ref{triplecorrelation}, one applies Proposition \ref{triplecorrelation} with 
\begin{align*}
S_1(f)&=S_2(f)=\widetilde{\Lambda}_R(f),\\
S_3(f)&=\widetilde{H}_R(f),
\end{align*}
and takes $X_1=R$ and $X_2=X_3=Q$, so that $\widehat{X}=\widehat{Q}=\widehat{N}^{1\slash 12}$,
\begin{align*}
S_{1,X_1}(f)&=\Lambda_R(f),\\
S_{2,X_2}(f)&=\Lambda_Q(f),\\
S_{3,X_3}(f)&=H_Q(f),
\end{align*}
and $y_1=1, y_2=-1$ and $y_3=0$. 

We have thus reduced our task to proving Proposition \ref{triplecorrelation}. To this end, we need the following lemma, along with some specialized estimates on sums involving the divisor function. 
\begin{lem}\label{lem:divisorspecial}
Let $m_1,m_2\in \F\setminus \{0\}$ be distinct. Let $r\in \F\setminus\{0\}$ and $b\in (\F\slash (r))^\times$. Let $h_1,h_2$ be square-free and suppose $g\in \F$ with $g\mid [h_1,h_2]r$. Then, 
\[\sum_{\substack{a_1\in (\F/(h_1))^\times\\a_2\in (\F/(h_2))^\times\\a_1h_2r+a_2h_1r\equiv bh_1h_2\hspace*{-2.5mm}\pmod{g(h_1,h_2)}}}e\left(\frac{m_1a_1}{h_1}+\frac{m_2a_2}{h_2}\right)\ll_{m_1,m_2}\tau(g)^2.\]
\end{lem}
\begin{proof}
Let us begin by setting $h_i':=h_i\slash (h_1,h_2)$ for $i=1,2$, $r'=r\slash (g,r)$, and $d=g\slash (g,r)$. If there are $a_1\in (\F/(h_1))^\times$ and $a_2\in (\F/(h_2))^\times$ satisfying the congruence  
\[a_1h_2r+a_2h_1r\equiv bh_1h_2\pmod{g(h_1,h_2)},\]
then we may divide through by $(h_1,h_2)$ to obtain 
\[a_1h_2'r+a_2h_1'r\equiv b w [h_1,h_2]\pmod{g},\]
where $w\in \mathbb{F}_q^{\times}$. This forces $(g,r)\mid b[h_1,h_2]$ so that $b[h_1,h_2]=v(g,r)$ for some $v\in\F$. Dividing the congruence above through by $(g,r)$ yields 
\[a_1h_2'r'+a_2h_1'r'\equiv v\pmod{d}.\]
Using the fact that $(r',d)=1$, we have 
\[a_1h_2'+a_2h_1'\equiv u\pmod{d},\]
where $u=v(r')^{-1}\in\F/(d)$. 

Now, our goal is to show 
\[S:=\sum_{\substack{a_1\in (\F/(h_1))^\times\\a_2\in (\F/(h_2))^\times\\a_1h_2'+a_2h_1'\equiv u\hspace*{-2.5mm}\pmod{d}}}e\left(\frac{m_1a_1}{h_1}+\frac{m_2a_2}{h_2}\right)\ll_{m_1,m_2}\tau(g)^2.\]
Using the orthogonality relation, we may rewrite $S$ as follows: 
\[
S=
\frac{1}{\langle d \rangle }\sum_{\lambda\in \F\slash (d)}e\bigg(\frac{-\lambda u}{d}\bigg) \sum_{\substack{a_1\in (\F/(h_1))^\times\\a_2\in (\F/(h_2))^\times}}e\left(\frac{m_1a_1}{h_1}+\frac{m_2a_2}{h_2}+\frac{\lambda(a_1h_2'+a_2h_1')}{d}\right).
\]
Note that $d\mid [h_1,h_2]$, so that $t:=[h_1,h_2]\slash d$ is in $\F$. It is easy to see that $h_2'\slash d=t\slash h_1$ and $h_1'\slash d=t\slash h_2$. Thus, 
\[S=\frac{1}{\langle d \rangle }\sum_{\lambda\in \F\slash (d)}e\bigg(\frac{-\lambda u}{d}\bigg)\sum_{\substack{a_1\in (\F/(h_1))^\times\\a_2\in (\F/(h_2))^\times}}e\left(\frac{a_1(m_1+\lambda t)}{h_1}\right)e\left(\frac{a_2(m_2+\lambda t)}{h_2}\right).\]
Using Corollary \ref{ramanujansumbound}, we obtain the bound 
\[|S|\leq \frac{1}{\langle d \rangle }\sum_{\lambda\in \F\slash (d)}\langle (h_1,m_1+\lambda t)\rangle \langle (h_2,m_2+\lambda t)\rangle.\]
If an irreducible polynomial $P$ divides $(h_1,m_1+\lambda t)$, then $P$ divides both $h_1$ and $m_1+\lambda t$. If $P$ divides $t$, then $P$ also divides $m_1$, and if $P$ does not divide $t$, then $P$ divides $h_1\slash (h_1,t)$. Thus, using the fact that $h_1$ is square-free, we see that 
\[(h_1,m_1+\lambda t)\mid m_1(h_1\slash (h_1,t), m_1+\lambda t).\]
In addition, it is not hard to see that $(h_1\slash (h_1,t))\mid d$. Therefore, 
\[m_1(h_1\slash (h_1,t), m_1+\lambda t)\mid m_1(d, m_1+\lambda t).\]
Putting these together, we deduce
\[\langle(h_1, m_1+\lambda t) \rangle\leq \langle m_1(d,m_1+\lambda t)\rangle. \]
Now, we have 
\[|S|\leq \frac{\langle m_1\rangle \langle m_2\rangle }{\langle d \rangle }\sum_{\lambda\in \F\slash (d)}\langle (d,m_1+\lambda t)\rangle \langle (d,m_2+\lambda t)\rangle.\]
Since $t=[h_1,h_2]\slash d$ and $h_1,h_2$ are square-free, we know that $(t,d)=1$. Thus,
\begin{align*}
&\frac{\langle m_1\rangle \langle m_2\rangle }{\langle d \rangle }\sum_{\lambda\in \F\slash (d)}\langle (d,m_1+\lambda t)\rangle \langle (d,m_2+\lambda t)\rangle\\
&=\frac{\langle m_1\rangle \langle m_2\rangle }{\langle d \rangle }\sum_{s\in \F\slash (d)}\langle (d,s)\rangle \langle (d,s+h)\rangle,
\end{align*}
where $h=m_2-m_1$ is nonzero by assumption. Let $d_1=(d,s)$ and $d_2=(d,s+h)$. Then $d_1,d_2\mid d$, and we have the congruences 
\begin{align*}
s&\equiv 0\pmod{d_1},\\
s&\equiv -h\pmod{d_2}.
\end{align*}
For fixed $d_1,d_2$, this system has at most one solution mod $[d_1,d_2]$ by the Chinese remainder theorem. Thus, there are at most $\langle d\slash[d_1,d_2]\rangle $ values of $s\pmod{d}$ which satisfy these congruences. Also, notice that $(d_1,d_2)\mid s$ and $(d_1,d_2)\mid s+h$ imply $(d_1,d_2)\mid h$, so that $\langle (d_1,d_2)\rangle \leq \langle h \rangle$. Piecing everything together, we get 
\begin{align*}
|S|\leq &\frac{\langle m_1\rangle \langle m_2\rangle }{\langle d \rangle }\sum_{s\in \F\slash (d)}\langle (d,s)\rangle \langle (d,s+h)\rangle\\
&=\frac{\langle m_1\rangle \langle m_2\rangle }{\langle d \rangle }\sum_{d_1,d_2\mid d}\langle d_1\rangle \langle d_2\rangle |\{s\in \F\slash (d): d_1=(d,s), d_2=(d,s+h)\}|\\
&\leq \frac{\langle m_1\rangle \langle m_2\rangle }{\langle d \rangle }\sum_{d_1,d_2\mid d}\langle d_1\rangle \langle d_2\rangle \left\langle \frac{d }{[d_1,d_2]}\right\rangle\\
&\leq \langle h \rangle \langle m_1\rangle \langle m_2\rangle \tau(d)^2.
\end{align*}
Since the last expression is $\ll_{m_1,m_2}\tau(g)^2$, we have the desired result. 
\end{proof}

Specifically, we need the following corollary. 
\begin{cor}\label{bound1}
Let $m_1,m_2,b,r\in \F\setminus \{0\}$ with $m_1,m_2$ distinct and $(b,r)=1$. Let $h_1,h_2$ be square-free. Then, uniformly in $h_3$,  
\[\sum_{\substack{a_1\in (\F/(h_1))^\times\\a_2\in (\F/(h_2))^\times\\\text{denom}(\frac{a_1}{h_1}+\frac{a_2}{h_2}-\frac{b}{r})=h_3}}e\left(\frac{m_1a_1}{h_1}+\frac{m_2a_2}{h_2}\right)\ll_{m_1,m_2}\tau(r[h_1,h_2])^3.\]
\end{cor}
\begin{proof}
The condition 
\[\text{denom}\left(\frac{a_1}{h_1}+\frac{a_2}{h_2}-\frac{b}{r}\right)=h_3\]
is equivalent to 
\[(a_1h_2r+a_2h_1r-bh_1h_2,h_1h_2r)=\frac{h_1h_2r}{h_3}.\]
(which, of course, requires $h_3\mid [h_1,h_2]r$) By Lemma \ref{Mobius} in Appendix, with \[f=a_1h_2r+a_2h_1r-bh_1h_2,\] $g=h_1h_2r$ and $\Delta=h_1h_2r\slash h_3$, we have 
\begin{align*}
&\1_{\text{denom}(\frac{a_1}{h_1}+\frac{a_2}{h_2}-\frac{b}{r})=h_3}\\
&=\1_{(a_1h_2r+a_2h_1r-bh_1h_2,h_1h_2r)=\frac{h_1h_2r}{h_3}}\\
&=\sum_{\substack{(h_1h_2r\slash h_3)\mid d\mid h_1h_2r\\ d~\text{monic}}}\mu\left(\frac{d}{h_1h_2r\slash h_3}\right)\1_{a_1h_2r+a_2h_1r\equiv bh_1h_2\hspace*{-2.5mm}\pmod{d}}\\
&=\sum_{\substack{([h_1,h_2]r\slash h_3)\mid d'\mid [h_1,h_2]r\\ d'~\text{monic}}}\mu\left(\frac{d'}{[h_1,h_2]r\slash h_3}\right)\1_{a_1h_2r+a_2h_1r\equiv bh_1h_2\hspace*{-2.5mm}\pmod{d'(h_1,h_2)}},
\end{align*}
where we have made the substitution $d'=d\slash (h_1,h_2)$. Now, we multiply this by $e\big(\frac{m_1a_1}{h_1}+\frac{m_2a_2}{h_2}\big)$ and sum over $a_1\in (\F/(h_1))^\times$ and $a_2\in (\F/(h_2))^\times$. By applying Lemma \ref{lem:divisorspecial} and the fact that there are at most $\tau(r[h_1,h_2])$ monic polynomials $d'$ dividing $r[h_1,h_2]$, we arrive at the desired result. 
\end{proof}

The proof of Proposition \ref{triplecorrelation} also relies on the following divisor sum estimate. 
\begin{lem}\label{bound2}
Let $B\geq 0$ and $r\in \F\setminus\{0\}$. Then given any $\epsilon>0$,
\[\sum_{\substack{\max\{\langle h_1\rangle, \langle h_2\rangle, \langle h_3\rangle\}\geq \widehat{X}\\ \mu(h_1)^2=\mu(h_2)^2=\mu(h_3)^2=1\\ h_1\mid rh_2h_3\\ h_2\mid rh_1h_3\\ h_3\mid rh_1h_2}}\frac{\tau(h_1)^B}{\langle h_1\rangle }\frac{\tau(h_2)^B}{\langle h_2\rangle }\frac{\tau(h_3)^B}{\langle h_3\rangle }\ll_{B,\epsilon} \langle r \rangle^{\epsilon}\widehat{X}^{\epsilon-1}.\]
\end{lem}
\begin{proof}
By symmetry, it suffices to consider only $h_1,h_2, h_3$ for which 
\[\langle h_1\rangle =\max\{\langle h_1\rangle ,\langle h_2\rangle, \langle h_3\rangle\}.\]
From the condition $h_1\mid rh_2h_3$, we have $\frac{h_1}{(h_1,rh_2)}\mid h_3$. In addition, 
\[(h_1,rh_2)=(h_1,r)\left(\frac{h_1}{(h_1,r)},h_2\right),\]
so
\[\langle h_3\rangle \geq \frac{\langle h_1\rangle }{\langle (h_1,r)\rangle \langle(\frac{h_1}{(h_1,r)},h_2) \rangle }.\]
And, there are at most $q\tau(r)\tau(h_1)\tau(h_2)$ polynomials $h_3$ which satisfy the condition $h_3\mid rh_1h_2$, and for any $h_3$ satisfying the condition we have $\tau(h_3)\leq \tau(r)\tau(h_1)\tau(h_2)$. Therefore,

\begin{align*}
&\sum_{\substack{\max\{\langle h_1\rangle, \langle h_2\rangle, \langle h_3\rangle\}\geq \widehat{X}\\ \mu(h_1)^2=\mu(h_2)^2=\mu(h_3)^2=1\\ h_1\mid rh_2h_3\\ h_2\mid rh_1h_3\\ h_3\mid rh_1h_2}}\frac{\tau(h_1)^B}{\langle h_1\rangle }\frac{\tau(h_2)^B}{\langle h_2\rangle }\frac{\tau(h_3)^B}{\langle h_3\rangle }\\
&\ll \tau(r)^{B+1} \sum_{\langle h_1\rangle \geq \widehat{X}}\frac{\tau(h_1)^{2B+1}\langle (h_1,r)\rangle}{\langle h_1\rangle^2 }\sum_{\langle h_2\rangle \leq \langle h_1\rangle }\frac{\tau(h_2)^{2B+1}\langle(\frac{h_1}{(h_1,r)},h_2) \rangle}{\langle h_2\rangle }
\\
&\ll_B  \tau(r)^{B+1} \sum_{\langle h_1\rangle \geq \widehat{X}}\frac{\tau(h_1)^{4B+3}\langle (h_1,r)\rangle \log^{O_B(1)}(\langle h_1\rangle)}{\langle h_1\rangle^2 },
\end{align*}
where in the last step we applied Lemma \ref{rankinbound} in Appendix to the inner sum. To finish off the last sum we apply Lemma \ref{rankinbound} once more on the dyadic intervals $\big[2^j\widehat{X},2^{j+1}\widehat{X}\big)$ for integers $j\ge0$. Applying Lemma \ref{rankinbound}, we have 
\begin{align*}
\sum_{\langle h_1\rangle\in[2^j\widehat{X},2^{j+1}\widehat{X})}\frac{\tau(h_1)^{4B+3}\langle (h_1,r)\rangle \log^{O_B(1)}(\langle h_1\rangle )}{\langle h_1\rangle^2 }&\leq \tau(r)^{4B+4}(2^j\widehat{X})^{-1}\log^{O_B(1)}(2^{j+1}\widehat{X}).
\end{align*}
Summing this on $j\ge0$ and applying the divisor bound $\tau(r)\ll_{\epsilon}\langle r\rangle^{\epsilon}$ yields the lemma. 
\end{proof}

We are finally able to prove Proposition \ref{triplecorrelation} by combining Corollary \ref{bound1} with Lemma \ref{bound2}.  

\begin{proof}[Proof of Proposition \ref{triplecorrelation}]
Let $\lambda=b\slash r$. We are trying to prove 
\begin{align}
\sum_{\max\big\{\frac{\langle h_1\rangle }{\widehat{X}_1},\frac{\langle h_2\rangle }{\widehat{X}_2},\frac{\langle h_3\rangle }{\widehat{X}_3}\big\}>1}&\alpha_1(h_1)\alpha_2(h_2)\alpha_3(h_3)\notag\\
&\times\sum_{\substack{a_1\in (\F/(h_1))^\times\\a_2\in (\F/(h_2))^\times\\a_3\in (\F/(h_3))^\times\\\frac{a_1}{h_1}+\frac{a_2}{h_2}+\frac{a_3}{h_3}\equiv\frac{b}{r}\hspace*{-2.5mm}\pmod{1}}}e\left(\frac{y_1a_1}{h_1}+\frac{y_2a_2}{h_2}+\frac{y_3a_3}{h_3}\right)\ll \langle r\rangle^{\epsilon}\widehat{X}^{\epsilon-1}\label{Goal}.
\end{align}
If 
\[\frac{a_1}{h_1}+\frac{a_2}{h_2}+\frac{a_3}{h_3}\equiv\frac{b}{r}\pmod{1},\]
then we clearly have $h_1\mid rh_2h_3, h_2\mid rh_1h_3,$ and $h_3\mid rh_1h_2$. In addition, given $a_1,a_2,h_1,h_2,h_3,b,r$, there is exactly one $a_3\in (\F/(h_3))^\times$ such that 
\[\frac{a_1}{h_1}+\frac{a_2}{h_2}+\frac{a_3}{h_3}\equiv\frac{b}{r}\pmod{1}\]
if and only if $\text{denom}(\frac{a_1}{h_1}+\frac{a_2}{h_2}-\frac{b}{r})=h_3$.
Therefore, we may rewrite the left-hand side of \eqref{Goal} as 
\begin{align*}
e\left(y_3\frac{b}{r}\right)\sum_{\substack{\max\big\{\frac{\langle h_1\rangle }{\widehat{X}_1},\frac{\langle h_2\rangle }{\widehat{X}_2},\frac{\langle h_3\rangle }{\widehat{X}_3}\big\}>1\\h_1\mid rh_2h_3\\ h_2\mid rh_1h_3\\h_3\mid rh_1h_2}}&\alpha_1(h_1)\alpha_2(h_2)\alpha_3(h_3)\\
&\times\sum_{\substack{a_1\in (\F/(h_1))^\times\\a_2\in (\F/(h_2))^\times\\\text{denom}(\frac{a_1}{h_1}+\frac{a_2}{h_2}-\frac{b}{r})=h_3}}e\left(\frac{(y_1-y_3)a_1}{h_1}+\frac{(y_2-y_3)a_2}{h_2}\right).
\end{align*}
By Corollary \ref{bound1} and the fact that $y_1,y_2,y_3$ are distinct, this is 
\[\ll\sum_{\substack{\max\big\{\frac{\langle h_1\rangle }{\widehat{X}_1},\frac{\langle h_2\rangle }{\widehat{X}_2},\frac{\langle h_3\rangle }{\widehat{X}_3}\big\}>1\\h_1\mid rh_2h_3\\ h_2\mid rh_1h_3\\h_3\mid rh_1h_2}}|\alpha_1(h_1)\alpha_2(h_2)\alpha_3(h_3)|\tau(r[h_1,h_2])^3.\]
Since $|\alpha_i(h_i)|\leq\tau(h_i)^B/\langle h_i\rangle $ and $\alpha_i$ is supported on square-free polynomials for $i=1,2,3$, this is 
\[\ll \tau(r)^3\sum_{\substack{\max\big\{\frac{\langle h_1\rangle }{\widehat{X}_1},\frac{\langle h_2\rangle }{\widehat{X}_2},\frac{\langle h_3\rangle }{\widehat{X}_3}\big\}>1\\\mu(h_1)^2=\mu(h_2)^2=\mu(h_3)^2=1\\h_1\mid rh_2h_3\\ h_2\mid rh_1h_3\\h_3\mid rh_1h_2}}\frac{\tau(h_1)^{B+3}}{\langle{h_1}\rangle}\frac{\tau(h_2)^{B+3}}{\langle{h_2}\rangle}\frac{\tau(h_3)^{B+3}}{\langle{h_3}\rangle}.\]
Applying Lemma \ref{bound2} and the divisor bound completes the proof. 
\end{proof}
\section{Ruzsa's construction}
We devote this section to a short proof of the lower bound in Theorem \ref{MainTheorem}, which adapts that of \cite[Theorem 4]{Ruzsa1984}. We consider the cases $\text{char}(\mathbb{F}_q)=2$ and $\text{char}(\mathbb{F}_q)>2$ separately.

\textbf{Case 1.} $\text{char}(\mathbb{F}_q)=2$. Fix an arbitrary $\epsilon>0$, and suppose $N\ge 2q^{3+\epsilon}$. Let $\{P_j\}_{j=1}^{J}\subseteq\F$ be the sequence of non-constant, monic irreducible polynomials of degree not exceeding
\[n=\left\lfloor\frac{\log((q-1)N/q)}{\log q}-\epsilon\right\rfloor\ge3,\] 
arranged in such a way that $\deg P_j\le\deg P_{j+1}$ for all $1\le j\le J-1$, where 
\[J=\sum_{k=1}^n\pi(k)=\sum_{k=1}^n\left(\frac{q^k}{k}+O\left(\frac{q^{k/2}}{k}\right)\right)=\frac{q^{n+1}}{(q-1)n}+O\left(\frac{q^{n}}{n^2}\right)\] 
by \eqref{eq:PNT}. We define a sequence $\{f_j\}_{j=1}^{J}\subseteq\F$ of monic polynomials inductively as follows. Let $f_1$ be the monic polynomial of least degree such that 
\[f_1\equiv -1\pmod{P_1}\quad\text{and}\quad f_1\equiv 0\pmod{P_2\cdots P_J}.\]
Once $f_{j-1}$ is determined, we take $f_{j}$ to be the monic polynomial of least degree, such that $\deg f_{j}>\deg f_{j-1}$, and
\[f_{j}\equiv -1\pmod{P_{j}}\quad\text{and}\quad f_{j}\equiv 0\pmod{P_{j+1}\cdots P_J}.\]
Consider the set 
\[A=\left\{g=\sum_{j=1}^{J}\epsilon_jf_j\colon\epsilon_1,...,\epsilon_J\in\{0,1\}\right\}.\]
Since $\deg f_{j}>\deg f_{j-1}$, we have $|A|=2^J$. By our construction of $\{f_j\}_{j=1}^{J}$, we know that 
\[\deg f_1\le\deg(P_1\cdots P_J)=\sum_{k=1}^n\pi(k)k=\sum_{k=1}^n\big(q^k+O\big(q^{k/2}\big)\big)=\frac{q^{n+1}}{q-1}+O\big(q^{n/2}\big),\]
and that 
\[\deg f_{j}\le \max\{1+\deg f_{j-1},\deg(P_1\cdots P_J)\}\]
for all $2\le j\le J$. It follows that 
\[\deg f_{j}\le\frac{q^{n+1}}{q-1}+O\big(q^{n/2}\big)+J\le N\]
for all $1\le j\le J$, provided that $N$ is sufficiently large in terms of $\epsilon$. Hence, every polynomial in $A$ has degree at most $N$ for sufficiently large $N$ in terms of $\epsilon$. Moreover, 
\[|A|=2^J\ge\exp\left((\log2+o(1))\frac{q^{n+1}}{(q-1)n}\right)\ge\exp\left((\log2+O(\epsilon))\frac{N\log q}{q\log N}\right)\]
for sufficiently large $N$ in terms of $\epsilon$. 

It suffices to show that the set $A$ constructed above has the desired property. Now, if $g,h\in A$ are distinct such that $h-g-f_k$ has degree less than $\deg f_k$ for some $1\le k\le J$, then 
\[h-g+1=(f_k+1)+\sum_{j=1}^{k-1}\delta_jf_j\equiv 0\pmod{P_k},\]
with each $\delta_j\in\{0,\pm1\}$, where we have observed that our assumption $\text{char}(\mathbb{F}_q)=2$ tells us that $f_k=-f_k$. Since $J\ge n\ge3$, we see that
\[\deg(h-g+1)=\deg(f_k+1)\ge\deg f_1\ge\deg(P_2\cdots P_{J})>\deg P_k.\]
It follows that $h-g+1$ is a nontrivial multiple of $P_k$ and hence reducible. We have thus shown that $A$ contains no pair of elements whose difference is of the form $P-1$ with $P$ irreducible. This completes the proof of the lower bound asserted in Theorem \ref{MainTheorem} when $\text{char}(F_q)=2$. \\

\textbf{Case 2.} $\text{char}(\mathbb{F}_q)>2$. The proof in this case is similar, so we will be brief. Again, let us fix an arbitrary $\epsilon>0$, and suppose $N\ge 2q^{6+\epsilon}$. Let $\{P_j\}_{j=1}^{J}\subseteq\F$ be the sequence of non-constant, monic irreducible polynomials of degree not exceeding $n$, where
\[n=\left\lfloor\frac{\log((q-1)N/q)}{2\log q}-\frac{\epsilon}{2}\right\rfloor\ge3,\] 
ordered in such a way that $\deg P_j\le\deg P_{j+1}$ for all $1\le j\le J-1$, where 
\[J=\frac{q^{n+1}}{(q-1)n}+O\left(\frac{q^{n}}{n^2}\right)\] 
as before. For $M=\lfloor J/2\rfloor$, we define a sequence $\{f_m\}_{m=1}^{M}\subseteq\F$ of monic polynomials inductively as follows. Let $f_1$ be the monic polynomial of least degree such that 
\[f_1\equiv -1\pmod{P_1}\quad\text{and}\quad f_1\equiv 1\pmod{P_2}.\]
Once $f_{m-1}$ is determined, we take $f_{m}$ to be the monic polynomial of least degree, satisfying the inequality $\deg f_{m}>\deg f_{m-1}$ and the congruences
\[f_{m}\equiv -1\pmod{P_{2m-1}},\quad f_{m}\equiv 1\pmod{P_{2m}},\quad\text{and}\quad f_{m}\equiv 0\pmod{P_{1}\cdots P_{2m-2}}.\]
The set
\[A:=\left\{g=\sum_{m=1}^{M}\epsilon_mf_m\colon\epsilon_1,...,\epsilon_M\in\{0,1\}\right\}\]
has cardinality
\[|A|=2^M\ge\exp\left(\left(\frac{\log2}{2}+O(\epsilon)\right)\frac{N\log q}{q\log N}\right)\]
for sufficiently large $N$ in terms of $\epsilon$. Moreover, we have
\[\deg f_{m}\le\max\{1+\deg f_{m-1},\deg(P_1\cdots P_{2m})\}\]
for all $2\le m\le M$, from which we deduce that
\[\deg f_{m}\le\frac{q^{2n+1}}{q-1}+O\big(q^{n}\big)+M\le N\]
for all $1\le m\le M$, assuming $N$ is sufficiently large in terms of $\epsilon$. So every polynomial in $A$ has degree at most $N$ for sufficiently large $N$ in terms of $\epsilon$. 

To verify that the set $A$ lacks shifted irreducible differences $P-1$, we take two arbitrary distinct $g,h\in A$ and suppose that $h-g-f_k$ has degree less than $\deg f_k$ for some $1\le k\le M$. Then $h-g+1$ is divisible by $P_{2k-1}$ or $P_{2k}$. Since
\[\deg(h-g+1)=\deg(f_k+1)>\deg P_{2k}\ge \deg P_{2k-1},\]
it follows that either $h-g+1$ is a nontrivial multiple of $P_{2k-1}$ or $P_{2k}$ and therefore reducible. This is sufficient for confirming the desired lower bound when $\text{char}(\mathbb{F}_q)>2$.

\section{Appendix}

We need the following elementary divisor sum bound. 
\begin{lem}\label{rankinbound}
For any $B>0$ and nonzero $h\in\F$,
\[\sum_{0<\langle f \rangle \leq \widehat{X}}\frac{\tau(f)^B\langle (f,h) \rangle }{\langle f\rangle }\ll_B \tau(h)^{B+1}\log^{2^B}\widehat{X}.\]
\end{lem}
\begin{proof}
Making use of the easy convolution identity
\[\sum_{\substack{d\mid g\\d\text{~monic}}}\phi(d)=\langle g \rangle\]
along with the sub-multiplicativity of $\tau$, we find that
\begin{align*}
\sum_{0<\langle f \rangle \leq \widehat{X}}\frac{\tau(f)^B\langle (f,h) \rangle }{\langle f\rangle}&=\sum_{\substack{d\mid h\\d\text{~monic}}}\phi(d)\sum_{\substack{0<\langle f \rangle \leq \widehat{X}\\d\mid f}}\frac{\tau(f)^B}{\langle f\rangle}\\
&\le(q-1)\sum_{\substack{d\mid h\\d\text{~monic}}}\frac{\phi(d)\tau(d)^B}{\langle d \rangle}\sum_{\substack{\langle g \rangle \leq \widehat{X}/\langle d\rangle\\g\text{~monic}}}\frac{\tau(g)^B}{\langle g\rangle},
\end{align*}
and by Euler factorization,
\begin{align*}
\sum_{\substack{\langle g \rangle \leq \widehat{X}/\langle d\rangle\\g\text{~monic}}}\frac{\tau(g)^B}{\langle g\rangle}&\le\prod_{\substack{\langle P\rangle
\leq \widehat{X}\slash \langle d \rangle \\P~\text{monic, irreducible}}}\left(1+\sum_{\nu\ge1}\frac{(\nu+1)^B}{\langle P \rangle^{\nu}}\right)\\
&= \prod_{\substack{\langle P \rangle \leq \widehat{X}\slash \langle d \rangle \\P\text{~monic, irreducible}}}\left(1+\frac{2^B}{\langle P\rangle}+O\left(\frac{1}{\langle P\rangle^2}\right)\right)\\
&=\exp\left(\sum_{\substack{\langle P\rangle \leq \widehat{X}\slash \langle d \rangle \\P\text{~monic, irreducible}}}\left(\frac{2^B}{\langle P\rangle}+O\left(\frac{1}{\langle P\rangle^2}\right)\right)\right).\\
\end{align*}
By the Prime Number Theorem for $\F$, we have that 
\begin{align*}
\sum_{\substack{\langle P\rangle \leq \widehat{X}\slash \langle d \rangle\\P\text{~monic, irreducible}}}\frac{1}{\langle P\rangle}=\sum_{q^n\le\widehat{X}\slash \langle d \rangle}\frac{1}{q^n}\sum_{\substack{\deg P=n\\P\text{~monic, irreducible}}}1&=\sum_{q^n\le\widehat{X}\slash \langle d \rangle}\left(\frac{1}{n}+O\left(\frac{q^{-n/2}}{n}\right)\right)\\
&=\log\log\big(\widehat{X}\slash \langle d \rangle\big)+ O(1),    
\end{align*}
and that
\[\sum_{\substack{\langle P\rangle \leq \widehat{X}\slash \langle d \rangle \\P\text{~monic, irreducible}}}\frac{1}{\langle P\rangle^2}\ll\sum_{q^n\le\widehat{X}\slash \langle d \rangle }\frac{1}{nq^n}\ll1.\]
But it is evident that
\[\sum_{\substack{d\mid h\\d\text{~monic}}}\frac{\phi(d)\tau(d)^B}{\langle d \rangle}\le\sum_{\substack{d\mid h\\d\text{~monic}}}\tau(h)^B=\tau(h)^{B+1}.\]
The proof of the lemma is now complete.
\end{proof}

We also need the following form of Möbius inversion.
\begin{lem}\label{Mobius}
Let $f\in \F$ and let $g,\Delta\in \F\setminus\{0\}$, with $\Delta$ monic. Then, 
\[\1_{(f,g)=\Delta}=\sum_{\substack{\Delta\mid d\mid g\\d~\text{monic}}}\mu(d\slash \Delta)\1_{d\mid f}.\]
\end{lem}
\begin{proof}
If $\Delta\nmid (f,g)$, then no term on the right-hand side contributes, so both sides are zero. Otherwise, we know that $\Delta\mid (f,g)$. Then, we may write $f=\Delta f'$, $g=\Delta g'$ and $d=\Delta d'$, and rewrite the right-hand side as 
\[\sum_{\substack{\Delta\mid d\mid g\\d~\text{monic}}}\mu(d\slash \Delta)\1_{d\mid f}=\sum_{\substack{d'\mid g'\\d'\text{~monic}}}\mu(d')\1_{d'\mid f'}=\sum_{\substack{d'\mid (f',g')\\d'\text{~monic}}}\mu(d')=\1_{(f',g')=1}=\1_{(f,g)=\Delta},\]
by the characteristic identity for the Möbius function. 

\end{proof}

\section{Acknowledgments} We would like to express our gratitude to Ákos Magyar, Giorgis Petridis, Paul Pollack, and Alex Rice for their valuable discussions and insights. We also thank the anonymous referees for a thorough and careful reading of the manuscript and for many helpful corrections and suggestions which led to substantial improvements in the clarity and presentation of our results.


\bibliographystyle{amsplain}
\nocite{*}
\bibliography{SarkozyFF}
\end{document}